\documentclass[twoside]{article}
\usepackage[accepted]{aistats2024}
\usepackage[round]{natbib}

\usepackage{amsmath}
\usepackage{amsfonts}
\usepackage{amssymb}
\usepackage{amsthm}
\usepackage{mathrsfs}
\usepackage{graphicx}
\usepackage{color}
\usepackage{hyperref}
\usepackage{bm}
\usepackage{algorithm}
\usepackage{algpseudocode}
\usepackage{algpascal}
\usepackage{url}
\usepackage{afterpage}
\usepackage{enumitem} 
\usepackage{multirow}
\usepackage{booktabs}
\usepackage{subfigure}
\usepackage{authblk}

\numberwithin{equation}{section}

\newtheorem{theorem}{Theorem}[section]
\newenvironment{Thm}{\begin{theorem}}{\end{theorem}}
\newtheorem{lemma}{Lemma}[section]
\newenvironment{Lem}{\begin{lemma}}{\end{lemma}}
\newtheorem{proposition}{Proposition}[section]

\newtheorem{corollary}{Corollary}[section]
\newenvironment{Cor}{\begin{corollary}}{\end{corollary}}
\theoremstyle{remark}

\newtheorem{example}{\bf Example}[section]

\renewenvironment{proof}{\noindent{\it Proof.}}{\qed}

\def\cC{\mathcal C}

\def\cF{\mathcal F}
\def\cG{\mathcal G}

\def\cN{\mathcal N}
\def\cP{\mathcal P}

\def\cZ{\mathcal Z}

\newcommand{\bb}{{\bf b}}

\newcommand{\bd}{{\bf d}}

\newcommand{\bff}{{\bf f}}
\newcommand{\bg}{{\bf g}}

\newcommand{\bh}{{\bf h}}

\newcommand{\bk}{{\bf k}}

\newcommand{\bfm}{{\bf m}}

\newcommand{\bt}{{\bf t}}

\newcommand{\bv}{{\bf v}}

\newcommand{\bx}{{\bf x}}
\newcommand{\bX}{{\bf X}}
\newcommand{\by}{{\bf y}}

\newcommand{\bz}{{\bf z}}
\newcommand{\bZ}{{\bf Z}}

\newcommand{\bbE}{{\mathbb E}}
\newcommand{\bbI}{{\mathbb I}}
\newcommand{\bbR}{{\mathbb R}}
\newcommand{\bbZ}{{\mathbb Z}}
\newcommand{\bbN}{{\mathbb N}}

\newcommand{\bbeta}{\bm{\beta}}

\newcommand{\bepsilon}{\bm{\epsilon}}

\newcommand{\iid}{{\rm i.i.d.}}

\newcommand{\Holder}{H\"{o}lder}

\newcommand{\bc}{\begin{center}}
	\newcommand{\ec}{\end{center}}
\newcommand{\be}{\begin{equation}}
	\newcommand{\ee}{\end{equation}}
\newcommand{\ba}{\begin{array}}
	\newcommand{\ea}{\end{array}}
\newcommand{\bean}{\setlength\arraycolsep{1pt}\begin{eqnarray*}}
	\newcommand{\eean}{\end{eqnarray*}}
\newcommand{\bea}{\setlength\arraycolsep{1pt}\begin{eqnarray}}
	\newcommand{\eea}{\end{eqnarray}}
\newcommand{\ben}{\begin{enumerate}}
	\newcommand{\een}{\end{enumerate}}
\newcommand{\bed}{\begin{itemize}}
	\newcommand{\eed}{\end{itemize}}

\DeclareMathOperator*{\argmax}{argmax}

\def\d{\mbox{d}}

\newcommand{\bzero}{{\bf 0}}

\newcommand{\vertiii}[1]{{\left\vert\kern-0.25ex\left\vert\kern-0.25ex\left\vert #1 
		\right\vert\kern-0.25ex\right\vert\kern-0.25ex\right\vert}}

\usepackage[symbol]{footmisc}

\begin{document}

	%
	\runningtitle{Minimax Optimal Density Estimation Using a Shallow Generative Model}
	
	%

	\twocolumn[
	
	\aistatstitle{Minimax Optimal Density Estimation Using a Shallow Generative Model with a One-Dimensional Latent Variable}
	
	\aistatsauthor{ Hyeok Kyu Kwon \And Minwoo Chae$^*$}
	\aistatsaddress{
		Department of Industrial and Management Engineering
		\\
		Pohang University of Science and Technology 
		\\
		$^*$Correspondence to: \href{mailto:mchae@postech.ac.kr}{mchae@postech.ac.kr} } ]

	\begin{abstract}
		A deep generative model yields an implicit estimator for the unknown distribution or density function of the observation. This paper investigates some statistical properties of the implicit density estimator pursued by VAE-type methods from a nonparametric density estimation framework. More specifically, we obtain convergence rates of the VAE-type density estimator under the assumption that the underlying true density function belongs to a locally \Holder\ class. Remarkably, a near minimax optimal rate with respect to the Hellinger metric can be achieved by the simplest network architecture, a shallow generative model with a one-dimensional latent variable. 
	\end{abstract}
	
	\section{INTRODUCTION}
	
	Suppose we have observations $\bX_1,\ldots, \bX_n$ that are $\iid$ copies of a $d$-dimensional random vector $\bX$ following the distribution $P_0$, with the density function $p_0$. Developing nonparametric estimators for $p_0$ has been a crucial task in unsupervised learning, and various methods and related theories are available in the literature \citep{hastie2009elements, tsybakov2008introduction, gine2016mathematical}. In recent years, deep generative models have shown remarkable success in modeling high-dimensional data, such as images and videos. Although classical density estimation methods provide direct estimators for $p_0$, deep generative model approaches can be seen as indirect estimation methods for $p_0$ because they only generate samples from the estimated distributions. Despite indirect estimation methods, deep generative models are very useful in many applications, including image and language generation problems.
	
	In our view, popularly used deep generative models can be categorized into two approaches based on their data-generating procedures. The first approach involves constructing an estimator $\hat\bg$ for a function $\bg: \bbR^{d_0} \to \bbR^d$, commonly referred to as the generator. Then, a sample $\bZ$ is drawn from a known $d_0$-dimensional distribution such as the standard normal or uniform, and $\hat\bg(\bZ)$ is treated as a sample from the estimated distribution. Thus, the distribution (or density) of $\hat\bg(\bZ)$ serves as an indirect estimator for $P_0$ (or $p_0$). Variational autoencoders (VAE) \citep{kingma2013auto, rezende2014stochastic}, normalizing flows (NF) \citep{dinh2015nice, rezende2015variational} and generative adversarial networks (GAN) \citep{goodfellow2014generative, arjovsky2017wasserstein, mroueh2018sobolev, li2017mmd} are important examples.
	
	The second approach involves estimating the score function, which is the gradient of the log density. Once an estimator of the score function is obtained, one can generate samples using score-based Markov chain Monte Carlo algorithms such as Hamiltonian and Langevin Monte Carlo \citep{neal2011mcmc}. Hence, the limit distribution of the Markov chain can be understood as an indirect estimator of $P_0$. The idea of score function estimation was originally suggested in \cite{hyvarinen2005estimation} and further developed in \cite{vincent2011connection, song2019generative, song2020sliced}. The score function estimation problem is closely related to the denoising diffusion model \citep{sohl2015deep, ho2020denoising}, and it has achieved state-of-the-art performance in many applications \citep{song2021scorebased}.
	
	Despite the tremendous success of deep generative models, their theoretical understanding remains largely unexplored. This paper focuses on studying the statistical theory for some generative model approaches. Specifically, we investigate the convergence rate of an implicit density estimator from a generative model. This estimator is the target estimator pursued by VAE approaches. Although it is empirically known in the literature that NF, GAN and score-based methods tend to outperform VAE, it deserves to study convergence rates of VAE type estimators because theoretical study provides a lot of valuable insights.
	
	Under the assumption that the true density $p_0$ belongs to a locally $\beta$-Hölder class, we prove that the estimator achieves the minimax optimal rate $n^{-\beta/(d + 2\beta)}$ up to a logarithmic factor with respect to the Hellinger metric. Remarkably, we show that the optimal rate can be achieved by the simplest ReLU \citep{glorot2011deep} network architecture consisting of a shallow network with a one-dimensional latent variable. Thus, even simple generative models can lead to optimal density estimators. The proof of the main theorem relies on the well-known result from the nonparametric Bayesian literature that a smooth density with a suitably decaying tail can efficiently be approximated by a finite mixture of normal distributions \citep{ghosal2001entropies, ghosal2007posterior, kruijer2010adaptive, shen2013adaptive}. The key is to find a tight upper bound for the number of support points of the mixing measure, which depends on the dimension and smoothness of the density.
	We also provide an alternative proof under additional assumptions, which offers important insights and suggests an extension to structured density estimation.
	This proof relies on the existence of a sufficiently regular generator for which Caffarelli's regularity theory of optimal transport \citep{caffarelli1990interior, villani2008optimal} provides sufficient conditions.
	
	There are several articles that investigate the convergence rates of implicit density estimators from deep generative models, with a focus on GAN-based approaches. \cite{liang2021well} and \cite{singh2018nonparametric} proved that a GAN-type estimator achieves the minimax optimal rate with respect to the Sobolev integral probability metric (IPM) \citep{muller1997integral}. The generalization to Besov IPMs can be found in \cite{uppal2019nonparametric}. \cite{belomestny2021rates} considered a vanilla GAN and obtained minimax optimal rates with respect to the Jensen--Shannon divergence. Note that all these results guarantee the optimal rate with respect to the total variation distance for sufficiently regular $p_0$. We would also like to mention earlier works \cite{pati2011posterior} and \cite{kundu2014latent}. Rather than parametrizing generators by neural networks, they considered Gaussian process priors and obtained optimal posterior convergence rates. Recently, diffusion models have also been considered in the context of implicit density estimation, and \cite{oko2023diffusion} obtained the minimax optimal rates with respect to the total variation and Wasserstein distances.
	
	Statistical theories for deep generative models beyond the nonparametric density estimation framework are also available in the literature, allowing for the possibility that $P_0$ is singular with respect to the Lebesgue measure. In this case, the parameter of interest is a distribution rather than a density. Various metrics have been considered to evaluate the performance of estimation, including the Sinkhorn divergence \citep{luise2020generalization}, Wasserstein metric \citep{chae2023likelihood, chae2022convergence} and general IPMs \citep{schreuder2021statistical, huang2021error, tang2023minimax, tang2024adaptivity}. These papers employ low-dimensional structures to explain how deep generative models can overcome the curse of dimensionality. For example, \cite{chae2023likelihood} and \cite{chae2022convergence} considered a composite structure on the generator, while \cite{tang2023minimax} assumed a manifold structure on the support of $P_0$ and derived the minimax optimal rate.

	The VAE-type estimator studied in this paper is analyzed in \cite{chae2023likelihood} under the assumption that $P_0$ is concentrated around a low-dimensional structure. Although the rate in \cite{chae2023likelihood} is not optimal, it is not significantly slower than the optimal rate, as discussed in \cite{chae2022convergence}. In contrast, the result in this paper guarantees that a VAE-type estimator is (nearly) optimal when $P_0$ has a smooth density. Combining these two results shows that, with carefully chosen network architectures, a VAE-type estimator can achieve a fast convergence rate regardless of the singularity of $P_0$. This highlights the adaptive nature of deep generative models to the structure of the unknown distribution.
	

	The remainder of this paper is organized as follows. In the following subsection, we provide notations and definitions. Section \ref{sec:model} introduces basic set-up and deep generative models. The main results concerning the convergence rate of VAE-type estimators are given in Section \ref{sec:main}.
	An alternative proof and extensions to the structured density estimation are given in Section \ref{sec:alternative}.
	Numerical results with a toy example and concluding remarks follow in Section \ref{sec:experiments} and \ref{sec:conclusion}, respectively.
	Technical proofs are provided in the supplementary material.

	\subsection{Notations and Definitions}
	
	A boldface is used to denote vectors. For $\bx \in \bbR^d$ and $1 \leq p \leq \infty$, let $\|\bx\|_{p}$ be the $\ell_{p}$-norm of $\bx$.
	For a set $A \subset \bbR^{d_1}$ and a vector-valued function $\bg = (g_1,\ldots,g_{d_2})^{\rm T} : A \rightarrow \bbR^{d_2}$, let
	\bean
	\|\bg\|_{p} = \left(\int_A \sum_{i=1}^{d_2} |g_{i}(\bz)|^{p} \d \bz \right)^{1/p}  \text{for $p \in [1, \infty)$},
	\eean
	and
	$
	\|\bg\|_{\infty} = \sup_{\bz \in A} \max(|g_1(\bz)|,\ldots,|g_{d_2}(\bz)|).
	$
	Let $\phi_{\sigma, d}$ be the density function of the multivariate normal distribution $\cN(\bm{0}_d, \sigma^2 \bbI_d )$, where $\bm{0}_d$ and $\bbI_d$ are $d$-dimensional zero vector and identity matrix, respectively. For simplicity, we often denote $\phi_{\sigma, d}$ as $\phi_\sigma$ when the dimension is obvious from the contexts.
	Let $\phi_{\sigma}*P$ be the convolution of $P$ and $\cN(\bzero_d, \sigma^2 \bbI_d)$, that is,
	\bean
	(\phi_{\sigma}*P)(\bx) = \int \phi_{\sigma}(\bx-\by) \d P(\by).
	\eean
	The Dirac measure at $\bx$ is denoted as $\delta_\bx(\cdot)$.
	For two probability density functions $p$ and $q$, the Kullback--Leibler (KL) divergence and Hellinger metric are denoted as
	\bean
	&& K(p, q) = \int p(\bx) \log \left(\frac{p(\bx)}{q(\bx)} \right) \d\bx
	\quad {\rm and} 
	\\
	&& d_H(p,q) = \left(\int \left\{ \sqrt{p(\bx)} - \sqrt{q(\bx}) \right\}^2 \d \bx \right)^{1/2},
	\eean
	respectively.
	For a (pseudo-)metric space $(\cP, \rho)$ and $\delta > 0$, let $N(\delta,\cP, \rho)$ and $N_{[]}(\delta,\cP, \rho)$ be the covering and bracketing numbers with respect to $\rho$, respectively. We refer to \cite{van1996weak} for details about these definitions. The notation $a \lesssim b$ implies that $a$ is less than or equal to $C b$, where $C$ is some constant that is not important in the given context. Similarly, $a \asymp b$ indicates that $a \lesssim b$ and $b \lesssim a$. Finally, the notation $C = C(A_1, \ldots, A_k)$ means that the constant $C$ depends solely on $A_1, \ldots, A_k$.

	\section{A LIKELIHOOD APPROACH TO DEEP GENERATIVE MODELS} \label{sec:model}
	
	This section presents a likelihood approach for deep generative models commonly used in practice. As previously mentioned, this method involves an estimator employed by VAE-type methods, which will henceforth be referred to as a VAE-type estimator.
	
	Our goal is to construct an estimator $\hat\bg$ of the generator $\bg: \bbR^{d_0} \to \bbR^d$ so that the distribution of $\hat\bg(\bZ)$ serves as an estimator of $P_0$, where $\bZ$ is a $d_0$-dimensional random vector following a known distribution. In particular, we aim to model $\bg$ using neural networks. Throughout this paper, we assume that $\bZ$ is a standard uniform variable on $[0,1]^{d_0}$. While likelihood-based approaches are a natural choice for constructing an estimator $\hat\bg$, deriving the likelihood for $\bg(\bZ)$ is difficult, and even the density of $\bg(\bZ)$ may not exist. Flow-based methods directly utilize the density of $\bg(\bZ)$, but this approach can limit the flexibility in designing network architectures.
	
	To overcome this difficulty, a VAE-type method employs an additional random vector and model $\bX$ as $\bX = \bg(\bZ) + \bepsilon$. Here, $\bepsilon$ is independent of $\bZ$ and follows the normal distribution $\cN(\bm{0}_d, \sigma^2 \bbI_d)$. Then, $\bX$ always allows the Lebesgue density 
	\be \label{eq:mixture-density}
	p_{\bg,\sigma}(\bx) = \int_{\left[0,1\right]^{d_0}} \phi_{\sigma}\left(\bx - \bg(\bz)\right) \d \bz
	\ee
	provided that $\sigma > 0$. Hence, one can obtain a maximum likelihood estimator by maximizing the log-likelihood function $(\bg, \sigma) \mapsto \sum_{i=1}^n \log p_{\bg, \sigma}(\bX_i)$ over $\cG \times [\sigma_{\min}, \sigma_{\max}]$, where $\cG$ is a class of functions from $[0,1]^{d_0}$ to $\bbR^d$ and $0 < \sigma_{\min} \leq \sigma_{\max} < \infty$. Formally, for a class $\cP$ of probability density functions and a sequence $(\eta_n)$ of nonnegative real numbers, an estimator $\hat p \in \cP$ is called an $\eta_n$-sieve MLE over $\cP$ if
	\bean
	\frac{1}{n} \sum_{i=1}^n \log \hat{p}(\bX_i) \geq \sup_{p \in \cP} \frac{1}{n} \sum_{i=1}^n \log p(\bX_i) - \eta_n.
	\eean
	Note that $\cP$, often called a sieve \citep{geman1982nonparametric}, is allowed to depend on the sample size, and $\eta_n$ can be understood as the optimization error. When $\cP$ consists of densities of the form \eqref{eq:mixture-density} with $\bg$ parametrized by deep neural networks, several algorithms approximating a sieve MLE have been suggested in the literature \citep{kingma2013auto, rezende2014stochastic, burda2016importance, dieng2019reweighted, kim2020casting}.
	
	To be more specific, for a positive integer $m$ and a vector $\bb = (b_1,\ldots,b_m)^{\rm T} \in \bbR^{m}$, let $\rho_{\bb}(\cdot) : \bbR^{m} \rightarrow \bbR^{m}$ be the ReLU activation function defined as
	\bean
	\rho_{\bb}(\bx) = \left(\max\{x_1-b_1,0\}, \ldots, \max\{x_m-b_m, 0\} \right)^{{\rm T}}
	\eean
	for $\bx = (x_1,\ldots,x_m)^{\rm T}$.
	For $L \in \bbN$, $F, M > 0$ and $\bd = (d_0,\ldots,d_{L+1}) \in \bbN^{L+2}$ with $d_{L+1}=d$, let $\cG = \cG(L,F,\bd,M)$ be the class of functions $\bg: [0,1]^{d_0} \to \bbR^d$ of the form
	\bean
	\bg(\bz) = W_{L}\rho_{\bb_{L}} \cdots W_1 \rho_{{\bb_1}} W_0 \bz  
	\eean
	with $W_i \in \bbR^{d_{i+1} \times d_{i}}$, $\bb_i \in \bbR^{d_i}$, $\|\bg\|_{\infty}\leq F$ and
	\bean
	\quad \max_{0 \leq i \leq L+1}\left\{ \max \left(\|W_i\|_{\infty} , \| \bb_i\|_{\infty} \right) \right\} \leq M,
	\eean
	where $\bb_0 = \bzero_{d_0}$ and $\|W_i\|_\infty$ is the entrywise maximum norm.
	
	In Section \ref{sec:main}, we analyze the convergence rate of an $\eta_n$-sieve MLE over
	\bean
	\cP = \Big\{p_{\bg,\sigma} : \bg \in \cG(L, F, \bd, M), \sigma \in [\sigma_{\rm min},\sigma_{\rm max}]  \Big\}
	\eean
	with $L=1$ and $\bd = (1, d_1, d)$. That is, the dimension of the latent variable $\bZ$ is 1, and the generator is parametrized by a shallow network with $d_1$ hidden units. Note that parameters such as $(F, d_1, M, \sigma_{\min})$ are allowed to depend on the sample size.

	\section{MAIN RESULTS} \label{sec:main}
	
	This section presents the main results of the paper. We first outline the assumptions on the true density $p_0$. Specifically, we will assume that $p_0$ belongs to a locally \Holder\ class with a suitably decaying tail. This class of density functions has been studied in \cite{shen2013adaptive} to analyze the convergence rate of the posterior distribution in a Dirichlet process mixture model. A slight improvement has been made in Chapter 9 of \cite{ghosal2017fundamentals}.
	
	\subsection{Assumptions on True Density Function} \label{ssec:assumption}
	
	For a multi-index $\bk = (k_1,\ldots,k_d)^{\rm T} \in (\bbZ_{\geq 0})^{d}$, denote $D^{\bk}$ the mixed partial derivative operator $\partial^{k.}/{\partial x_1^{k_1}\cdots \partial x_d^{k_d}}$, where $k. = \sum_{j=1}^d k_j$. 
	For any $\beta>0, \tau_0 \geq 0$ and non-negative function $L : \bbR^d \rightarrow \bbR$, let $\cC^{\beta,L,\tau_0} (A)$ be the class of every real-valued function $f$ on $A \subseteq \bbR^d$ such that
	$\sup_{\bx \in A} | D^{{\bk}} f(\bx) | < \infty$ for $ k. \leq \lfloor \beta \rfloor$, and 
	\bean
	|(D^{\bk}f)(\bx+\by) - (D^{\bk}f)(\bx)| \leq L(\bx) e^{\tau_0 \| \by \|_2^2} \|\by\|_2^{\beta - \lfloor \beta \rfloor}
	\eean
	for $k. = \lfloor \beta \rfloor,\bx \in A$ and $\by \in \{\bz:\bx+\bz \in A\}$, where $\lfloor \beta \rfloor$ denotes the largest integer strictly smaller than $\beta$.
	
	We will assume that $p_0 \in \cC^{\beta,L,\tau_0}(\bbR^d)$ for some $\beta$, $\tau_0$ and $L$. We also make the following two technical assumptions on the tail of $p_0$.
	
	\begin{itemize}[label={}, labelsep=0pt, leftmargin=0pt]
		\item (Tail 1) For any $\bk \in (\bbZ_{\geq 0})^{d}$ with $k. \leq \lfloor \beta \rfloor$,
		\bean
		\bbE \left[ \left( \frac{L(\bX)}{p_0(\bX)} \right)^{2} + \left( \frac{|D^{\bk} p_0(\bX)|}{p_0(\bX)} \right)^{\frac{2\beta}{k.}} \right] < \infty,
		\eean
		where $\bbE$ denotes the expectation with respect to $P_0$.
		\item (Tail 2) There exist $\tau_1,\tau_2, \tau_3 > 0$ such that $p_0(\bx) \leq \tau_1  \exp(-\tau_2 \|\bx\|_2^{\tau_3})$ for all $\bx \in \bbR^d$.
	\end{itemize}
	
	The above assumptions, in particular the tail assumptions, are satisfied by a large class of densities. For example, suppose that $p_0$ is the $d$-dimensional standard normal density. Then, for any $\bk \in (\bbZ_{\geq 0})^d$, we have $D^\bk p_0(\bx) \lesssim (1 + \|\bx\|_1)^{k.} p_0(\bx)$ because the standard normal density $\phi$ satisfies $\phi'(x) = -x \phi(x)$. Therefore, 
	\bean
	&& |(D^\bk p_0) (\bx + \by) - (D^\bk p_0)(\bx) | 
	\\
	&& = \left| \by^T \int_0^1 \nabla (D^\bk p_0) (\bx + t\by) \d t \right| 
	\\
	&& \lesssim \|\by \|_2 \sup_{t \in [0,1]} \Big[ (1 + \|\bx + t\by\|_1)^{k.+1} p_0(\bx + t\by) \Big]
	\\
	&& \lesssim
	\|\by \|_2 (1+\|\bx \|_1 + \| \by\|_1)^{k.+1} e^{ \frac{\|\by\|_2^2}{2\alpha} - \frac{\|\bx\|_2^2}{2(1+\alpha)}}
	\eean
	for every $\alpha > 0$, where the last inequality holds because $\|\bx + t\by\|_2^2 \geq \|\bx\|_2^2/(1+\alpha) - \|\by\|_2^2/\alpha$ for all $t \in [0,1]$. Hence, for any $\beta > 0$, if we take $\alpha = 1/2$, $\tau_0 > 1$ and $L(\bx) = c (\|\bx\|_1^{\lfloor \beta\rfloor +1} + 1) e^{-\|\bx\|_2^2/3}$ for a large enough constant $c = c(\beta, d, \tau_0)$, then $p_0 \in \cC^{\beta,L,\tau_0}(\bbR^d)$ and two tail conditions are satisfied with $\tau_1 = (2\pi)^{-d/2}$, $\tau_2 = 1/2$ and $\tau_3 = 2$.
	
	As another example, suppose that $p_0$ is the $d$-fold product density of the Laplace distribution, that is, $p_0(\bx) = 2^{-d}e^{-\|\bx\|_1}$. Simple calculation yields that 
	\bean
	&& \left\vert p_0(\bx+\by) - p_0(\bx) \right\vert = 2^{-d}e^{-\|\bx\|_1} \left\vert 1-e^{-\|\bx+\by\|_1 + \|\bx\|_1} \right\vert 
	\\
	&& \leq 2^{-d} \|\by\|_1  e^{-\|\bx\|_1}
	\eean
	for $\bx,\by \in \bbR^d$, where the inequality holds because $1 - e^{-x} \leq x$ for all $x \in \bbR$.
	Since $\|\by\|_1 \leq \sqrt{d}\|\by\|_
	2$, $p_0$ belongs to $\cC^{1, L,0}(\bbR^{d})$ with $L(\bx) = \sqrt{d}2^{-d}e^{-\|\bx\|_1}$. Furthermore, two tail conditions are satisfied with $\tau_1 = 2^{-d}, \tau_2 = 1$ and $\tau_3 =1$ because $\|\bx\|_2 \leq \|\bx\|_1$.
	\subsection{Convergence Rate of a Sieve MLE}
	\label{ssec:convergence}
	
	Under the assumptions stated in Section \ref{ssec:assumption}, it has been proven in \cite{shen2013adaptive} (and Chapter 9 of \cite{ghosal2017fundamentals}) that the posterior distribution, which is based on the Dirichlet location mixture of normal prior with a Gaussian base measure and an inverse Wishart prior on the covariance matrix parameter, contracts to $p_0$ with a minimax rate up to a logarithmic factor. An important technique used is to approximate $p_0$ by a finite mixture of normal distributions. The following lemma summarizes the result, and its proof can be easily derived from Lemmas 9.11 and 9.12 of \cite{ghosal2017fundamentals}. Hereafter, $C = C(\rm all)$ means that $C$ is a constant depending only on $d, \beta, L$ and $\tau_j$'s.
	
	\medskip
	
	\begin{lemma} \label{lem:informal}
		For any density function $p_0 \in \cC^{\beta, L, \tau_0}(\bbR^{d})$ satisfying assumptions (Tail 1) and (Tail 2), and small enough $\sigma > 0$, there exists a discrete probability measure $H(\cdot) = \sum_{i=1}^{N} w^{(i)} \delta_{\bx^{(i)}}(\cdot)$ supported within a compact set $E_{\sigma} = [-C \{\log (1/\sigma)\}^{\tau_3}, C \{\log (1/\sigma)\}^{\tau_3}]^d$ such that
		\bean 
		d_H(p_0,\phi_{\sigma} * H) \lesssim \sigma^{\beta} \left\{ \log (1/\sigma) \right\}^{d/4}
		\eean
		and $N \lesssim \sigma^{-d} \{ \log (1/\sigma)\}^{\tau_3 d+ d}$, where $C = C({\rm all})$.
	\end{lemma} 
	
	The approximation error improves as the smoothness of the density $p_0$ increases, according to Lemma \ref{lem:informal}. This lemma has been used in \cite{shen2013adaptive} to construct a sieve with metric entropy suitably bounded. We utilize it to approximate $p_0$ by a density of the form \eqref{eq:mixture-density} with $\bg$ a shallow ReLU network. Theorem \ref{thm:main} below is our main result.

	\medskip
	
	\begin{Thm} \label{thm:main}
		Suppose that $p_0 \in \cC^{\beta, L, \tau_0}(\bbR^d)$ and assumptions (Tail 1) and (Tail 2) are satisfied.
		Then, there exists a constant $\widetilde C_0 = \widetilde C_0({\rm all})$ such that for every constant $\widetilde C \geq \widetilde C_0$, an $\eta_n$-sieve MLE $\hat p$ over
		\bean
		\cP = \Big\{p_{\bg,\sigma} : \bg \in \cG(1,F,\bd,M), \sigma \in [\sigma_{\rm min},\sigma_{\rm max}] \Big\},
		\eean
		with $\bd = (1,d_1,d)$, $\sigma_{\min} = n^{-1/(2\beta+d)}$, $\sigma_{\max} = 1$ and
		\bean
		&& F = \widetilde C \left( \log n \right)^{\tau_3}, 
		\ d_1 = \big\lfloor \widetilde C n^{\frac{d}{2\beta+d} } \left( \log n \right)^{\tau_3 d+d} \big\rfloor,
		\\
		&& M = \widetilde{C} n^{\frac{2\beta+2d+3}{2\beta+d}},
		\eean
		satisfies 
		\bean
		P_0 \Big(d_H(p_0,\hat{p}) > \epsilon_n \Big) \leq 5 \exp \left(- A n \epsilon_n^2 \right) +  n^{-1} \log n
		\eean
		for every $n \geq \widetilde C_1$, where $\widetilde C_1=\widetilde C_1({\rm all}, \widetilde C)$, $\widetilde C_2=\widetilde C_2({\rm all}, \widetilde C)$, $\eta_n = \epsilon_n^2 / 48$,
		\bean
		\epsilon_n = \widetilde{C}_2 n^{-\frac{\beta}{2\beta+d}}  \left( \log n \right)^{\frac{2\tau_3 d+2\tau_3 +2d + 1}{2}}
		\eean
		and $A> 0$ is an absolute constant.
	\end{Thm}
	
	The statement of Theorem \ref{thm:main} has strong restrictions on the model parameters due to our attempt to minimize unimportant constants. However, it can be inferred from the proof that the parameters can be chosen more flexibly. For instance, one can choose $\sigma_{\min}=n^{-c_1}$ for a constant $c_1>1/(2\beta+d)$, $\sigma_{\max}=c_2$ for a constant $c_2\geq1$, $F = n^{c_3}$ for a constant $c_3 > 0$, and $M=n^{c_4}$ for a constant $c_4 > (2\beta+2d+3)/(2\beta+d)$. The key is to control the order of $d_1$, which determines the approximation and estimation errors for the density estimation.

	The proof of Theorem \ref{thm:main} involves several technical details and is provided in the supplementary material. Here, we provide an overview of the key ideas behind the proof. For convenience, we use the informal notation $a \lesssim_{\log} b$ to indicate that $a$ is less than or equal to $b$ up to a poly-logarithmic factor, such as $\log n$, $|\log \sigma|^d$, and $|\log \epsilon_n|^{\tau_3}$. Similarly, we use the notation $\asymp_{\log}$.
	
	To establish a convergence rate for the sieve MLE over the class $\cP$, we rely on the general theory developed in \cite{wong1995probability}, specifically Theorem 4. In essence, Theorem 4 states that a sieve MLE can achieve a suitable convergence rate if the KL divergence between the true density $p_0$ and the class $\cP$ is small enough and the bracket entropy of $\cP$ is suitably bounded. More specifically, if
	\be \begin{split} \label{eq:wongshen}
		& \inf_{p \in \cP} K(p_0, p) \lesssim_{\log} \epsilon_n^2
		\quad {\rm and}
		\\
		& \log N_{[]} (\epsilon_n, \cP, d_H) \lesssim_{\log} n\epsilon_n^2,
		\end {split} \ee
		then a sieve MLE over $\cP$ attains a convergence rate of $\epsilon_n$ with respect to the Hellinger metric. Note that each inequality is used to bound the approximation and estimation errors. Since $K(p_0, p) \asymp_{\log} d_H^2 (p_0, p)$ under a mild integrability condition (see Theorem 5 of \cite{wong1995probability} and Lemma B.2 of \cite{ghosal2017fundamentals}), the first inequality in \eqref{eq:wongshen} can be replaced by $\inf_{p \in \cP} d_H(p_0, p) \lesssim_{\log} \epsilon_n$.
		
		If we take $\sigma \asymp n^{-1/(2\beta+d)}$ in Lemma \ref{lem:informal}, we have
		\bean
		d_H(p_0,\phi_{\sigma}*H) \lesssim_{\log} \epsilon_n
		\quad {\rm and} \quad
		N \lesssim_{\log} n\epsilon_n^2,
		\eean
		where $H(\cdot) = \sum_{i=1}^{N} w^{(i)} \delta_{\bx^{(i)}}(\cdot)$ is the discrete measure in Lemma \ref{lem:informal}. Therefore, it suffices to show that the density function $\phi_\sigma * H$ can be approximated by the class $\cP$ of shallow ReLU network functions, with an approximation error of $\epsilon_n$ with respect to the Hellinger metric and bracket entropy of $n\epsilon_n^2$. For this purpose, we construct a ReLU network function $\bg: [0,1] \to \bbR^d$ so that the distribution of $\bg(Z)$ is sufficiently close to the discrete measure $H$, where $Z$ is a standard uniform random variable. 
		
\begin{figure*}[!t]
			\centering
			\subfigure[$\phi_\sigma*H = p_{\widetilde{\bg}, \sigma}$ for some $\widetilde{\bg}$]{
				\includegraphics[width=0.45\textwidth]{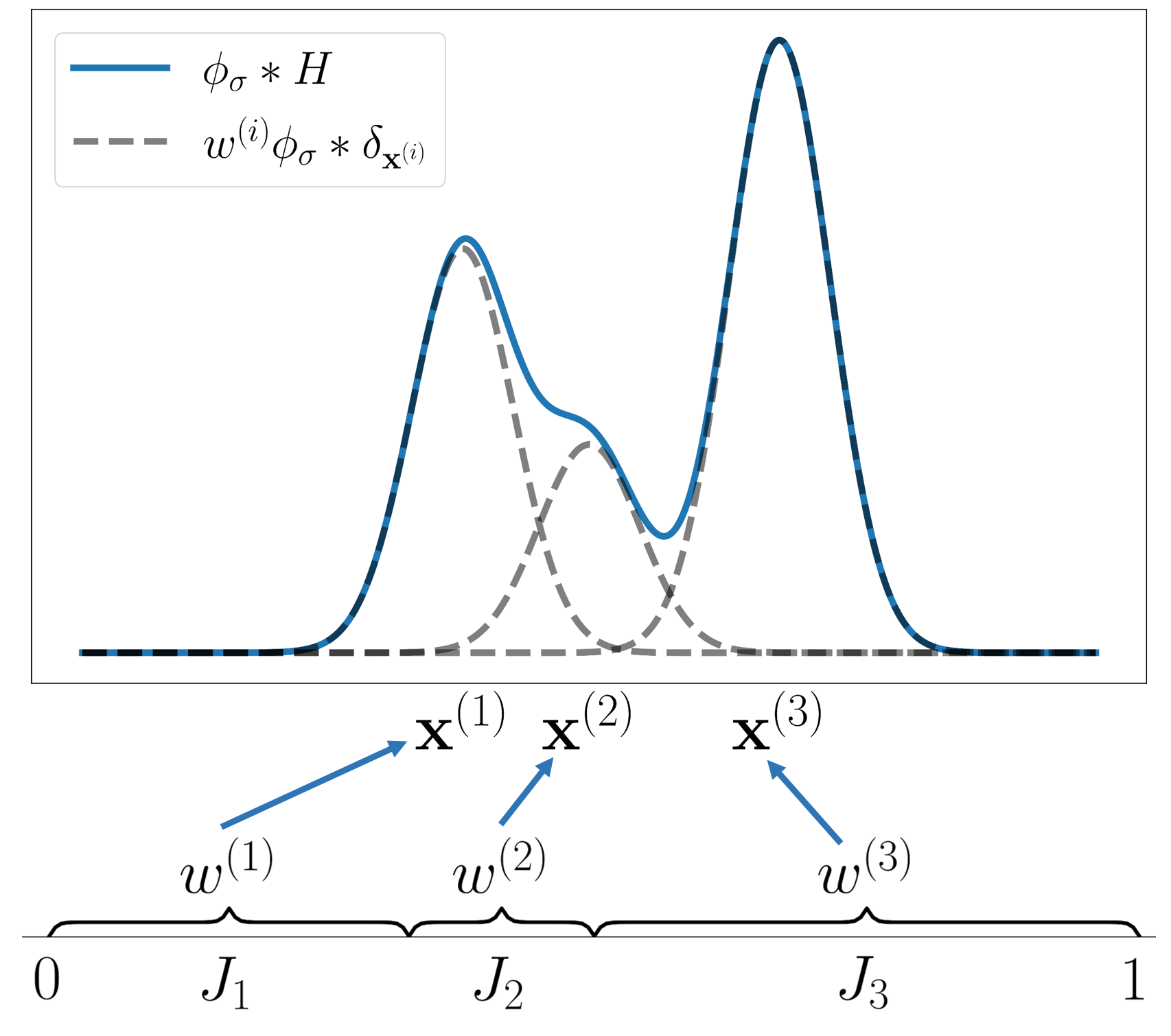}
				\label{Fig:gaussianmixture}
			}
			\subfigure[Approximation of $\widetilde{\bg}$ by a ReLU network $\bg$]{
				\includegraphics[width=0.51\textwidth]{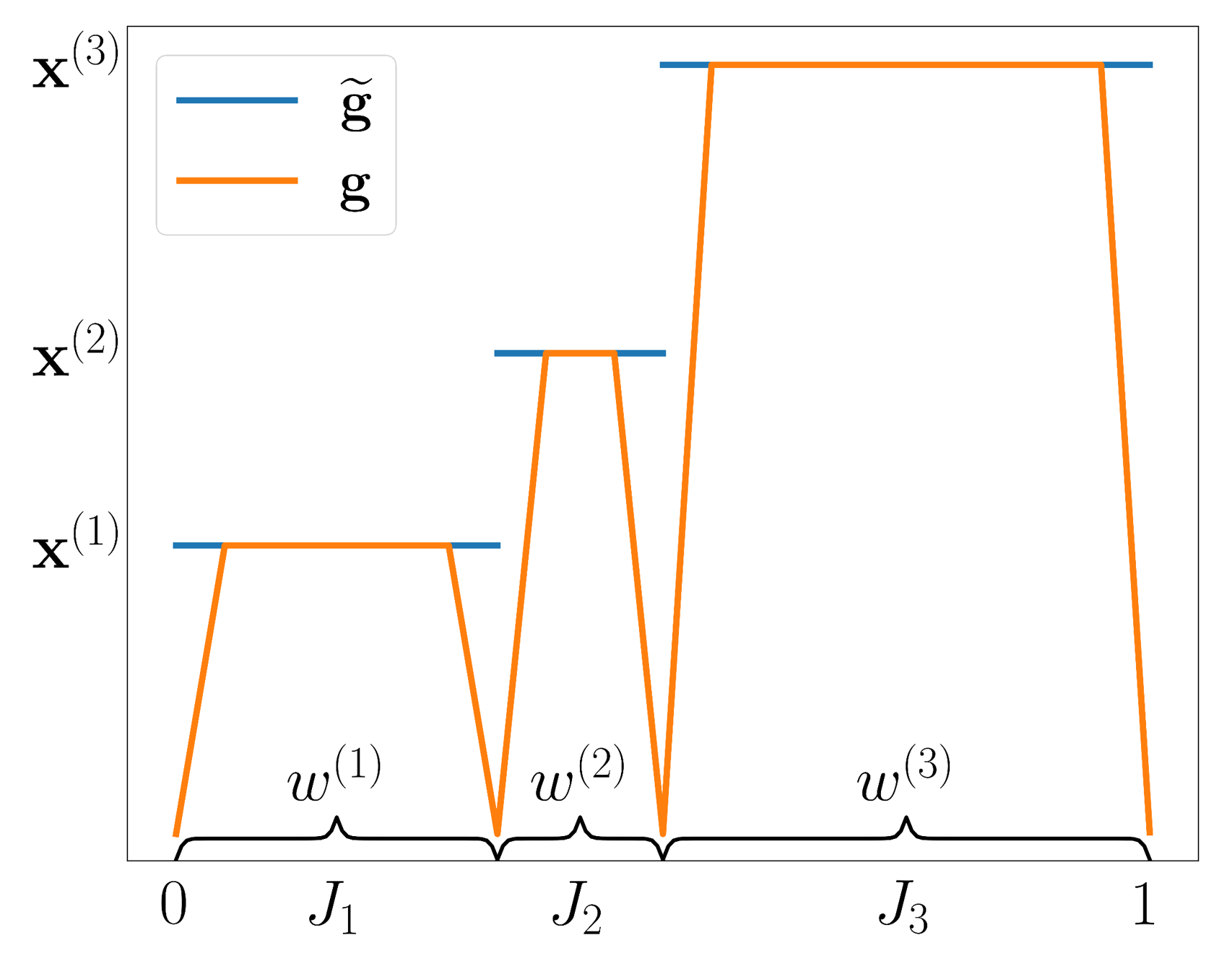}
				\label{Fig:relunetworks}
			}
			
			\caption{
				(a) A finite mixture $\phi_\sigma * H$ can be represented as $P_{\widetilde{\bg}, \sigma}$ for some function $\widetilde{\bg}: [0,1] \to \bbR^d$.
				(b) If $\widetilde{\bg}$ is a sum of $N$ indicator functions, it can be approximated by a shallow ReLU network function with $O(N)$ units.
			} \label{fig:approx}
		\end{figure*}

		The main idea of constructing such a $\bg$ is illustrated in Figure \ref{fig:approx}. We first define $\widetilde{\bg} (z) = \sum_{i=1}^N \bx^{(i)} 1_{J_i}(z)$ for consecutive intervals $J_1, \ldots, J_N$ that partition the unit interval $[0,1]$, where $\mu(J_i) = w^{(i)}$, and $\mu$ denotes the Lebesgue measure. It is easy to see that $H$ equals the distribution of $\widetilde{\bg}(Z)$. Next, we approximate each summand $\bx^{(i)} 1_{J_i}(\cdot)$, which is a constant function on the interval $J_i$, with a piecewise linear function, or equivalently, a shallow ReLU network. Since $\widetilde{\bg}$ is the sum of $N$ indicator functions, the number of hidden units required for the shallow ReLU approximation is of order $O(N)$. Therefore, by defining $\cP$ as in Theorem \ref{thm:main}, we can achieve the first inequality of \eqref{eq:wongshen}. Since the number $d_1$ of hidden units is of order $O(N) \lesssim_{\log} O(n\epsilon_n^2)$, the log of the $\epsilon_n$-covering number of the shallow network class $\cG(1, F, \bd, M)$ with respect to the uniform norm $\|\cdot\|_\infty$ is also of order $O(N)$ up to a logarithmic factor. This leads to the bracket entropy bound in \eqref{eq:wongshen}, completing the proof of Theorem \ref{thm:main}.
		
		It is worth noting that while Theorem \ref{thm:main} is limited to the ReLU activation function, other choices of activation functions are possible. From the previous sketch of the proof, we can see that the primary role of neural networks is to approximate the indicator function $1_{J_i}(\cdot)$. As ReLU networks are piecewise linear, they can easily approximate $1_{J_i}(\cdot)$ as in Figure \ref{fig:approx}-(b). Although not as straightforward as the ReLU activation function, it is possible for other activation functions to approximate $1_{J_i}(\cdot)$. In particular, Lemma 4 of \cite{imaizumi2022adaptive} shows that commonly used activation functions such as Sigmoid, LeakyReLU \citep{maas2013rectifier}, SoftPlus \citep{dugas2000incorporating}, and Swish \citep{ramachandran2017searching} can also approximate $1_{J_i}(\cdot)$ well. Thus, these activation functions can replace the ReLU in Theorem \ref{thm:main}.
		
		An important implication of Theorem \ref{thm:main} is that the minimax optimal rate for nonparametric density estimation can be achieved by the simplest network architecture. While some mathematical properties of shallow networks have been studied in the literature, most of them focus on the approximation properties of network functions. It is well-known that shallow networks have universal approximation capability \citep{cybenko1989approximation}. Furthermore, \cite{mhaskar1996neural} obtained nearly optimal numbers of hidden units to approximate a smooth function with the sigmoidal activation function. Although \cite{mhaskar1996neural} did not consider a statistical problem, the approximation theory might lead to optimal convergence rates for statistical problems, such as nonparametric regression. Recently, \cite{yang2023optimal} proved that shallow networks with ReLU activation function can lead to an optimal rate in nonparametric regression. To the best of our knowledge, the mathematical and statistical properties of shallow generative models, particularly those with a one-dimensional latent variable, have not been studied in the literature.

		\section{AN ALTERNATIVE PROOF AND STRUCTURED DENSITY ESTIMATION} 
		\label{sec:alternative}
		
		In this section, we present an alternative approach to obtain the convergence rate in Theorem \ref{thm:main} using deep generative models rather than shallow networks. While this alternative approach requires additional assumptions, it sheds light on potential extensions to structured density estimation and provides valuable insights. Moreover, our investigation has revealed a potential limitation of VAE-type estimators for structured density estimation.
		

		In addition to the $\beta$-regularity of $p_0$, we assume the existence of a $(\beta+1)$-regular function $\bg_0: \cZ \to \bbR^d$ such that $P_0$ is the distribution of $\bg_0(\bZ)$, where $\bZ$ follows a known distribution supported on $\cZ \subset \bbR^d$. The regularity theory of optimal transport by \cite{caffarelli1990interior} provides a sufficient condition for the existence of $(\beta+1)$-regular $\bg_0$ under the assumption that $p_0$ is $\beta$-regular. See Theorem 12.50 of \cite{villani2008optimal} for a general and rigorous statement, and \cite{cordero2019regularity} for state-of-the-art results. These statements involve several intricate notions from the Monge-Amp\`{e}re equation, so we also refer to Lemma 10 of \cite{chae2023likelihood} for readers who are not familiar with these notions. Note that the existence of a $(\beta+1)$-regular $\bg_0$ has been assumed in \cite{belomestny2021rates} to prove that the vanilla GAN achieves the minimax rate $n^{-\beta/(2\beta+d)}$ for nonparametric density estimation.

		For $\tau > 0$, the global \Holder\ class $\cC^{\beta}(A;\tau)$ is defined as the class of function $f \in \cC^{\beta,\tau,0}(A)$ satisfying $\sup_{\bx \in A} |D^{\bk} f(\bx)| \leq \tau$ for $k. \leq \lfloor \beta \rfloor$. For a vector valued function, we denote $\bff \in \cC^{\beta}(A;\tau)$ if each component of $\bff$ belongs to $\cC^{\beta}(A;\tau)$. Now, we specify additional assumptions used for the alternative approach. 
		\begin{itemize}[label={}, labelsep=0pt, leftmargin=0pt]
			\item (Support) There exists a constant $\tau_4 > 0$ such that $\{\bx: p_0(\bx) > 0  \} \subset [-\tau_4,\tau_4]^{d}$.
			\item (Generator) There exists a constant $\tau_5 \geq 1$ such that $P_0$ is the distribution of $\bg_0(\bZ)$ for some $\bg_0 \in \cC^{\beta+1}([0,1]^{d};\tau_5)$, where $\bZ$ is a uniform random vector on $[0,1]^d$.
		\end{itemize}
		
		We will also assume that $\beta \leq 2$ for technical reasons described below.
		Although it is unclear whether it is possible to achieve the minimax rate to the case $\beta > 2$, the case $\beta \leq 2$ is sufficient to discuss the benefit of the alternative approach and structured density estimation.
		Under these additional assumptions, we consider a sieve MLE $\hat p$ over $\cP = \{p_{\bg, \sigma}: \bg \in \cG\}$, where $\cG = \cG(L,F,\bd,M,s)$, the set of functions $\bg \in \cG(L, F, \bd, M)$ with the number of nonzero network parameters bounded by $s$.
		
		\begin{Thm} \label{thm:generator}
			Suppose that $p_0 \in \cC^{\beta, L,\tau_0}(\bbR^{d})$ with $\beta \leq 2$ and assumptions (Tail 1), (Support) and (Generator) are satisfied. Then, there exists a constant $\widetilde{C}_0 = \widetilde{C}_0(d,\beta,\tau_4,\tau_5)$ such that for every constant $\widetilde{C} \geq \widetilde{C}_0$, an $\eta_n$-sieve MLE $\hat{p}$ over 
			\bean
			\cP = \Big\{p_{\bg,\sigma}: \bg \in \cG(L,F,\bd,M,s), \sigma \in [\sigma_{\min},\sigma_{\max}] \Big\}
			\eean
			with $\bd = (d,d_1,\ldots,d_1,d) \in \bbN^{L+2}, \sigma_{\min} = n^{-1/(2\beta+d)}, \sigma_{\max} = 1$ and
			\bean
			&& L = \lfloor \widetilde{C} \log n \rfloor, \ F = \widetilde{C}, \ d_1 = \lfloor \widetilde{C} n^{\frac{d}{2\beta + d}} \rfloor, \\
			&& M = 1, \ s = \widetilde{C} n^{\frac{d}{2\beta + d}} \log n,
			\eean
			satisfies 
			\bean
			P_0 \Big(d_H(p_0,\hat{p}) > \epsilon_n \Big) \leq 5 \exp \left(- A n \epsilon_n^2 \right) +  n^{-1} 
			\eean
			for every $n \geq \widetilde C_1$, where $\widetilde C_1=\widetilde C_1({\rm all}, \widetilde C)$, $\widetilde C_2=\widetilde C_2({\rm all}, \widetilde C)$, $\eta_n = \epsilon_n^2 / 48$,
			\bean
			\epsilon_n = \widetilde{C}_2 n^{-\frac{\beta}{2\beta+d}} \log n
			\eean
			and $A> 0$ is an absolute constant.
		\end{Thm}
		
		Theorem \ref{thm:generator} is a special case of Theorem \ref{thm:structured}. Here, we only provide an overview of the key ideas behind the proof.
		By the well-known approximation property of deep neural networks \citep{schmidt2020nonparametric, ohn2019smooth, yarotsky2017error, telgarsky2016benefits}, there exists a network function $\bg \in \cG$ such that $\|\bg_0 - \bg \|_\infty \lesssim_{\log} s^{-(\beta+1)/d}$. Combining this with a convolution approximation $d_H(p_0, \phi_\sigma * P_0) \lesssim \sigma^\beta$ (see Lemma \ref{sec:beta2 approximation} and Chapter 4 of \cite{gine2016mathematical}) and Lemma \ref{sec:mixture hellinger} in the supplementary materials leads to an approximation error bound 
		\bean
		\inf_{p \in \cP} d_H(p_0, p) \lesssim_{\log} \sigma^\beta + \frac{s^{-(\beta+1)/d}}{\sigma}.
		\eean
		The $\delta$-entropy of $\cG$ with respect to the uniform metric is of order $O(s \log (1/\delta))$ up to a logarithmic factor, which provides a similar bound on the bracketing entropy of $\cP$. By choosing $s \asymp_{\log} n^{d/ (2\beta + d)}$ and $\sigma \asymp n^{-1/(2\beta+d)}$, the general approach of \cite{wong1995probability}, see also \eqref{eq:wongshen}, leads to the Hellinger convergence rate of $\epsilon_n \asymp_{\log} n^{-\beta/(2\beta+d)}$.
		
		Note that $d_H(p_0, \phi_\sigma * P_0) \lesssim \sigma^{\beta}$ does not hold for $\beta > 2$. For an extension to $\beta > 2$, more technical details should be involved as in \cite{kruijer2010adaptive} and \cite{shen2013adaptive}. We leave this as future work.
		
		Although the alternative approach requires additional assumptions, it can be used to develop a statistical theory that explains the benefits of deep generative models compared to shallow ones. Specifically, we consider structured density estimation, where the structure of a density is imposed through the generator. We assume that in addition to the regularity assumptions on $p_0$ and $\bg_0$, $\bg_0$ has a composite structure of the form
		\be \label{eq:composite}
		\bg_0 = \bh_{q} \circ \bh_{q-1} \circ \cdots \circ \bh_{1} \circ \bh_{0}
		\ee
		with $\bh_i = (h_{i1},\ldots,h_{iv_{i+1}})^{\rm T} : [a_i,b_i]^{v_i} \rightarrow [a_{i+1},b_{i+1}]^{v_{i+1}}$. Here, $v_0=v_{q+1}=d$ and $t_i$ is the maximal number of variables on which each component of $\bh_{i}$ depends. For any $q \in \bbZ_{\geq 0}, \bv = (v_0,\ldots,v_{q+1})^{\rm T} \in \bbN^{q+2}, \bt = (t_0,\ldots,t_{q})^{\rm T} \in \bbN^{q+1}, \bbeta = (\beta_0,\ldots,\beta_{q})^{\rm T} \in (\bbR_{> 0})^{q+1}$ and $  \tau > 0$, let $\cF(q,\bv,\bt,\bbeta,K)$ be the class of every real-valued functions of the form \eqref{eq:composite} satisfying $h_{ij} \in \cC^{\beta_i}([a_i,b_i]^{t_i};\tau)$ and $\max(|a_i|,|b_i|) \leq \tau$. Let
		\bean
		&& i_* = \argmax_{i \in \{0, \ldots,q\}} \frac{t_i}{\beta_i}, \quad \beta_* = \beta_{i_*}
		\quad
		\text{and} ~~ t_* = t_{i_*}.
		\eean
		Then, the assumption can be represented as follows.
		\begin{itemize}[label={}, labelsep=0pt, leftmargin=0pt]
			\item (Structured generator) $P_0$ is the distribution of $\bg_0(\bZ)$ for some  $\bg_0 \in \cF(q,\bv,\bt,\bbeta,\tau_6)$, with $\min_{i} \beta_i > 1$, where $\bZ$ is a uniform random vector on $[0,1]^d$.
		\end{itemize}
		This composite structure has been previously studied in the context of nonparametric regression by \cite{schmidt2020nonparametric} and \cite{bauer2019deep} to explain the benefits of deep neural networks. In the context of deep generative models, \cite{chae2023likelihood} and \cite{chae2022convergence} have used this structure to impose a low-dimensional structure on singular distribution estimation problems. 
		
		Similarly, we consider this composite structure on the generator for nonparametric structured density estimation. The general approach of \cite{wong1995probability} can still be used to obtain a convergence rate.
		\begin{Thm} \label{thm:structured}
			Suppose that $p_0 \in \cC^{\beta,L,\tau_0}(\bbR^{d})$ and assumptions (Tail 1), (Support) and (Structured generator) are satisfied. Let $\widetilde{\beta} = \min(\beta,2)$. Then, there exists a constant $\widetilde{C}_0 = \widetilde{C}_0(\beta,\tau_4,q,\bv,\bt,\bbeta,\tau_6)$ such that for every constant $\widetilde{C} \geq \widetilde{C}_0$, an $\eta_n$-sieve MLE $\hat{p}$ over 
			\bean
			\cP = \Big\{p_{\bg,\sigma}: \bg \in \cG(L,F,\bd,M,s), \sigma \in [\sigma_{\min},\sigma_{\max}] \Big\}
			\eean
			with $\bd = (d,d_1,\ldots,d_1,d) \in \bbN^{L+2}, \sigma_{\min} = n^{-\frac{\beta_*}{t_*(\widetilde{\beta}+1)+2\widetilde{\beta} \beta_*}},\sigma_{\max} = 1$ and
			\bean
			&& L = \lfloor \widetilde{C} \log n \rfloor, \ F = \widetilde{C}, \ d_1 = \lfloor \widetilde{C} n^{\frac{t_*(\widetilde{\beta}+1)}{2\widetilde{\beta}\beta_* + t_*(\widetilde{\beta}+1)}} \rfloor, \\
			&& M = 1, \ s = \widetilde{C} n^{\frac{t_*(\widetilde{\beta}+1)}{2\widetilde{\beta} \beta_* + t_*(\widetilde{\beta}+1)}} (\log n),
			\eean
			satisfies 
			\bean
			P_0 \Big(d_H(p_0,\hat{p}) > \epsilon_n \Big) \leq 5 \exp \left(- A n \epsilon_n^2 \right) +  n^{-1} 
			\eean
			for every $n \geq \widetilde C_1$, where $\widetilde C_1=\widetilde C_1({\rm all},q,\bv,\bt,\bbeta, \widetilde C)$, $\widetilde C_2=\widetilde C_2({\rm all},q,\bv,\bt,\bbeta, \widetilde C)$, $\eta_n = \epsilon_n^2 / 48$,
			\bean
			\epsilon_n = \widetilde{C}_2 n^{-\frac{\widetilde{\beta}\beta_*}{2\widetilde{\beta}\beta_* + t_*(\widetilde{\beta}+1)}}  \left( \log n \right)
			\eean
			and $A> 0$ is an absolute constant.
		\end{Thm}
		Note that Theorem \ref{thm:generator} is a special case of Theorem \ref{thm:structured} with $q=0,t_*=d$, $\beta_*=\beta+1$ and $\tau_5 = \tau_6$. Roughly speaking, a class $\cG$ of deep neural networks with $s$ nonzero parameters can approximate $\bg_0$ with an approximation error of $s^{-\beta_*/t_*}$. Also, the $\delta$-bracket entropy of $\cP = \{p_{\bg, \sigma}: \bg \in \cG\}$ can be bounded by $s \log(1/\delta)$ up to a logarithmic factor on $\sigma$. Hence, the general approach leads to the convergence rate $\epsilon_n \asymp_{\log} \sigma^{\widetilde{\beta}} + s^{-\beta_* / t_*}/\sigma + \sqrt{s/n}$. By taking
		\bean
		s \asymp n^{ \frac{t_* (\widetilde{\beta}+1)}{t_*(\widetilde{\beta}+1) + 2\widetilde{\beta}\beta_*} }
		\quad {\rm and} \quad
		\sigma \asymp s^{- \frac{\beta_*}{t_* (\widetilde{\beta}+1)}},
		\eean
		we obtain the Hellinger rate of
		\bean
		\epsilon_n \asymp_{\log} n^{-\frac{\widetilde{\beta}\beta_*}{2\widetilde{\beta}\beta_* + t_*(\widetilde{\beta}+1)}}.
		\eean
		Note that the rate depends on the dimension only through $t_*$, which might be much smaller than $d$. Additionally, it depends on both $\beta$ and $\beta_*$, where $\beta$ represents the smoothness of $p_0$ and $\beta_*$ is the smoothness of the worst component functions of $\bg_0$.
		
		
		The structured density estimation described above has not been studied in the literature; thus, the minimax optimal rate is unknown. It is worth noting that since we only need to estimate the generator $\bg_0$, it seems undesirable for the convergence rate to depend on $\beta$, the smoothness of $p_0$. However, with a VAE-type estimator considered in the present paper, the dependence on $\beta$ appears to be inevitable due to the convolution approximation error $d_H(p_0, \phi_\sigma * P_0) \lesssim \sigma^\beta$. The NF approach could be a promising alternative for obtaining the optimal rate because it directly utilizes the density of $\bg(\bZ)$. It is empirically known that NF outperforms VAE in many applications; therefore, in the future, it will be worth studying the convergence rate of NF approaches in structured density estimation.

		\section{NUMERICAL EXPERIMENTS}
		\label{sec:experiments}
		
		\begin{figure*}[!t]
			\centering
			\subfigure[ $d_H^2(\hat p, p_0)$ for KDE and VAE-type methods ]{
				\includegraphics[width=0.48\textwidth]{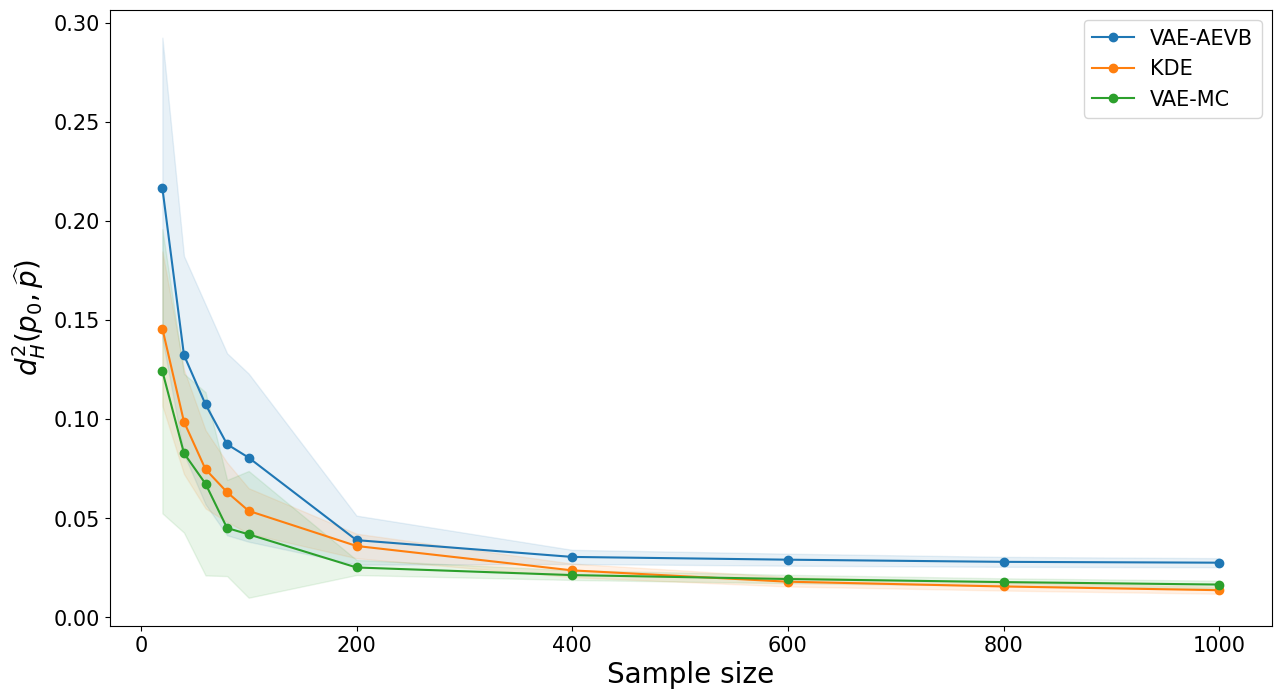}
			}
			\subfigure[Training log-likelihood values]{
				\includegraphics[width=0.48\textwidth]{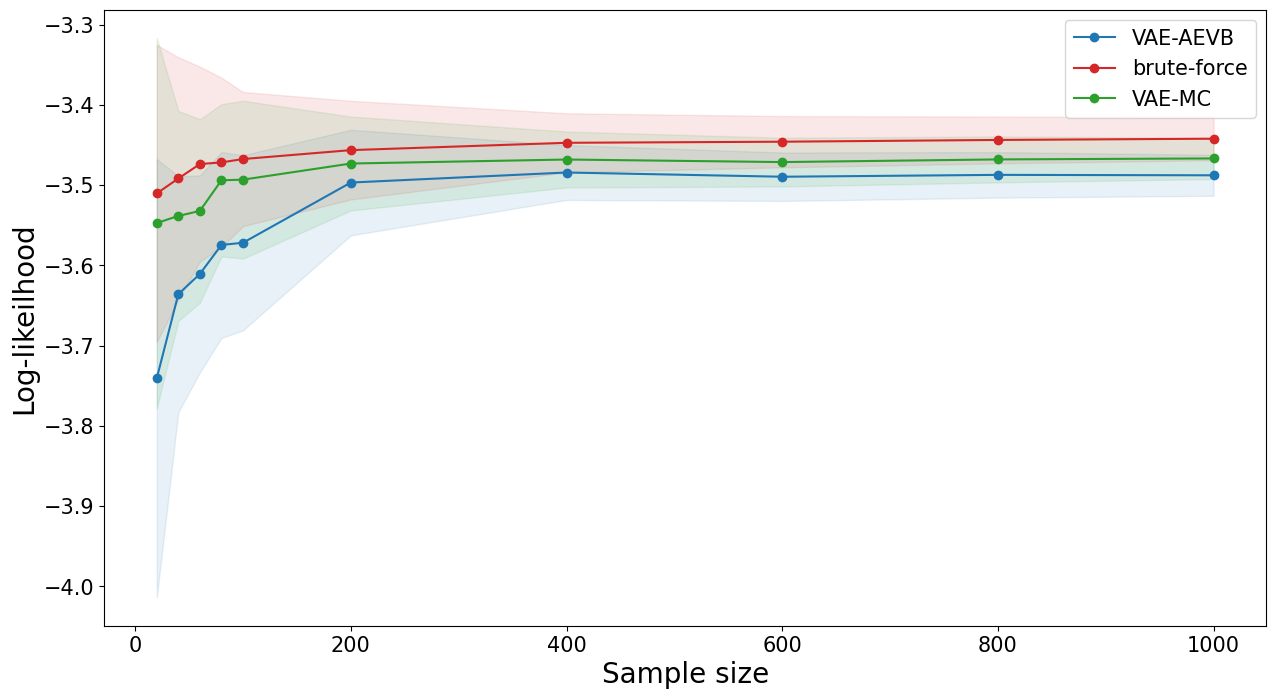}
				
			}
			
			\caption{The means and standard deviations of the squared Hellinger distances and training log-likelihood values. All results are based on 50 repetitions.}\label{fig:sim}
		\end{figure*}

		In this section, we conduct small numerical experiments to assess the actual performance of a shallow generative model with a one-dimensional latent variable.
		Data are generated from a two-component Gaussian mixture with $d=2$.
		More specifically, the true density is defined as $p_0(\cdot) = 0.5\phi(\cdot - \bfm) + 0.5\phi(\cdot + \bfm)$ with $\bfm = (1.3,1.3)^{\rm T}$.
		We consider a shallow ReLU network function $\bg_{\theta}$ parameterized by $\theta$ with $50$ hidden units.
		Since the likelihood function $p_{\bg_{\theta},\sigma}$ of the form (\ref{eq:mixture-density}) is computationally intractable, we approximate it using two approaches, Monte Carlo integration and auto-encoding variational Bayes (AEVB) algorithm \citep{kingma2013auto, rezende2014stochastic}. 
		
		For the Monte Carlo method, the log-likelihood is approximated as
		\bean
		\hat L_{\rm MC}(\theta,\sigma;\bx) = \log \left( \frac{1}{m} \sum_{i=1}^{m} \phi_{\sigma}(\bx-\bg_{\theta}(Z_i)) \right),
		\eean
		where $Z_1,\ldots,Z_{m}$ are standard uniform random variables.
		Then, one can obtain an implicit estimator $\hat{p}$ by maximizing $\sum_{i=1}^n \hat{L}_{\rm MC}(\theta,\sigma;\bX_i)$, which will be referred to as VAE-MC.

		Alternatively, one can maximize a lower bound of the log-likelihood using variational methods \citep{jordan1999introduction}. 
		Define the variational density $z \mapsto q_{\psi}(z \vert \bx)$ as the density of $\cN(\mu_{\psi}(\bx),\sigma_{\psi}^2(\bx))$, where $\mu_{\psi}(\cdot)$ and $2 \log \sigma_{\psi}(\cdot)$ are parameterized by neural networks, specifically as shallow ReLU networks with $50$ hidden units for the experiments. 
		For each iteration, define
		\bean
		\hat L_{{\rm AEVB}, i}(\theta,\sigma,\psi ; \bX_i) = \log \left( \frac{p_{\theta,\sigma}(\bX_i, Z_i)}{q_{\psi}(Z_i \vert \bX_i)} \right),
		\eean
		where $Z_i$ is a sample from $q_{\psi}(\cdot \vert \bX_i)$, $p_{\theta,\sigma}(\bx,z) = \phi_{\sigma}(\bx  -\bg_{\theta}(\Phi(z))) \phi(z)$ and $\Phi$ is the cumulative distribution function of the standard Gaussian distribution. 
		Then, one can obtain $\hat p$ by maximizing $\sum_{i=1}^n \hat L_{{\rm AEVB}, i}(\theta,\sigma,\psi; \bX_i)$, which will be referred to as VAE-AEVB.
		
		Both $\hat L_{\rm MC}$ and $\hat L_{\rm AEVB}$ are maximized using the Adam optimization algorithm \citep{kingma2014adam} with a mini-batch of size 20. The learning rate is fixed at $2 \times 10^{-4}$ for 1000 epochs and $m = 10^5$ is used for $\hat L_{\rm MC}$. 
		
		To evaluate the estimation performance, the squared Hellinger distance $d_H^2(\hat p, p_0)$ is computed for VAE-MC, VAE-AEVB and the Gaussian kernel density estimator (KDE). Silvermann's method is used to estimate the bandwidth parameter in KDE, implemented in Scikit-learn \citep{pedregosa2011scikit}.
		Note that the numerical integration implemented in SciPy \citep{virtanen2020scipy} is applied to compute the Hellinger distances.
		The results with varying sample sizes are depicted in Figure \ref{fig:sim}-(a).
		While VAE-MC performs comparably to KDE, VAE-AEVB performs significantly worse than KDE.
		This discrepancy is mainly due to the highly non-convex nature of the objective function used in the AEVB algorithm, leading to the failure of the SGD-based algorithm to maximize the log-likelihood. (Note that the number of parameters in the VAE-AEVB objective is about twice as great as in the VAE-MC objective.)
		To confirm this, we obtained network parameters with a high likelihood value using a brute-force method, which relies on the unknown structure of $p_0$. 
		The brute-force method sets $\sigma = 1$ and defines a piecewise linear function $\bg_{\theta}(\cdot)$ that closely approximates the sum of two indicator functions $\bfm \{1_{[0,0.5)}(\cdot) - 1_{[0.5,1]}(\cdot)\}$, as shown in Figure \ref{fig:approx}-(b). Specifically, $\bg_{\theta}$ is constructed as shallow ReLU networks with 8 hidden units as in (\ref{eq:nn_aprox}), with $\kappa = 10^{-5}$. Figure \ref{fig:sim}-(b) compares the training log-likelihood values of the VAE-type methods and brute-force method, confirming the failure of the SGD-based algorithm in maximizing the log-likelihood value.
		

		\section{CONCLUSIONS} \label{sec:conclusion}
		
		The VAE is an important class of inferential methods for deep generative models, but it is widely known that other methods, such as GAN, NF, and score-based methods, often outperform VAE in various applications. However, our paper shows that even the VAE with the simplest network architecture can produce a nearly optimal estimator in the nonparametric density estimation framework. This finding highlights the importance of considering further structures of the density or distribution being estimated to explain the superior performance of deep generative models over classical nonparametric methods.
		
		We suggest that the composite structure on the generator, as discussed in Section \ref{sec:alternative}, could be a promising structural assumption to investigate in future studies of density or distribution estimation problems. Such studies could lead to a better understanding of the benefits of deep generative models over classical nonparametric methods, and potentially inspire the development of even more powerful and efficient deep generative models.

		\section*{ACKNOWLEDGEMENT}
		
		The authors are grateful to the Area Chair and four anonymous reviewers for their valuable comments and suggestions on an earlier version of the paper. This work was supported by Samsung Science and Technology Foundation under Project Number SSTF-BA2101-03.
		
		\clearpage
		\bibliography{bib-short}

\begin{thebibliography}{}

\bibitem[Arjovsky et~al., 2017]{arjovsky2017wasserstein}
Arjovsky, M., Chintala, S., and Bottou, L. (2017).
\newblock Wasserstein generative adversarial networks.
\newblock In {\em Proc. International Conference on Machine Learning}, pages
  214--223.

\bibitem[Bauer and Kohler, 2019]{bauer2019deep}
Bauer, B. and Kohler, M. (2019).
\newblock On deep learning as a remedy for the curse of dimensionality in
  nonparametric regression.
\newblock {\em Ann. Statist.}, 47(4):2261--2285.

\bibitem[Belomestny et~al., 2021]{belomestny2021rates}
Belomestny, D., Moulines, E., Naumov, A., Puchkin, N., and Samsonov, S. (2021).
\newblock {Rates of convergence for density estimation with GANs}.
\newblock {\em ArXiv:2102.00199}.

\bibitem[Burda et~al., 2016]{burda2016importance}
Burda, Y., Grosse, R., and Salakhutdinov, R. (2016).
\newblock Importance weighted autoencoders.
\newblock In {\em Proc. International Conference on Learning Representations},
  pages 1--14.

\bibitem[Caffarelli, 1990]{caffarelli1990interior}
Caffarelli, L.~A. (1990).
\newblock {Interior $W^{2,p}$ estimates for solutions of the Monge--Amp\`{e}re
  equation}.
\newblock {\em Ann. of Math.}, 131(1):135--150.

\bibitem[Chae, 2022]{chae2022convergence}
Chae, M. (2022).
\newblock Rates of convergence for nonparametric estimation of singular
  distributions using generative adversarial networks.
\newblock {\em ArXiv:2202.02890}.

\bibitem[Chae et~al., 2023]{chae2023likelihood}
Chae, M., Kim, D., Kim, Y., and Lin, L. (2023).
\newblock A likelihood approach to nonparametric estimation of a singular
  distribution using deep generative models.
\newblock {\em J. Mach. Learn. Res.}, 24:1--42.

\bibitem[Cordero-Erausquin and Figalli, 2019]{cordero2019regularity}
Cordero-Erausquin, D. and Figalli, A. (2019).
\newblock Regularity of monotone transport maps between unbounded domains.
\newblock {\em Discrete and Continuous Dynamical Systems}, 39(12):7101--7112.

\bibitem[Cybenko, 1989]{cybenko1989approximation}
Cybenko, G. (1989).
\newblock Approximation by superpositions of a sigmoidal function.
\newblock {\em Math. Control Signals Systems}, 2(4):303--314.

\bibitem[Dieng and Paisley, 2019]{dieng2019reweighted}
Dieng, A.~B. and Paisley, J. (2019).
\newblock Reweighted expectation maximization.
\newblock {\em ArXiv:1906.05850}.

\bibitem[Dinh et~al., 2015]{dinh2015nice}
Dinh, L., Krueger, D., and Bengio, Y. (2015).
\newblock Nice: Non-linear independent components estimation.
\newblock In {\em Proc. International Conference on Learning Representations},
  pages 1--13.

\bibitem[Dugas et~al., 2000]{dugas2000incorporating}
Dugas, C., Bengio, Y., B\'{e}lisle, F., Nadeau, C., and Garcia, R. (2000).
\newblock Incorporating second-order functional knowledge for better option
  pricing.
\newblock In {\em Proc. Neural Information Processing Systems}, volume~13,
  pages 1--7.

\bibitem[Geman and Hwang, 1982]{geman1982nonparametric}
Geman, S. and Hwang, C.-R. (1982).
\newblock Nonparametric maximum likelihood estimation by the method of sieves.
\newblock {\em Ann. Statist.}, 10(2):401--414.

\bibitem[Ghosal and van~der Vaart, 2017]{ghosal2017fundamentals}
Ghosal, S. and van~der Vaart, A. (2017).
\newblock {\em Fundamentals of Nonparametric Bayesian Inference}.
\newblock Cambridge University Press.

\bibitem[Ghosal and van~der Vaart, 2001]{ghosal2001entropies}
Ghosal, S. and van~der Vaart, A.~W. (2001).
\newblock {Entropies and rates of convergence for maximum likelihood and Bayes
  estimation for mixtures of normal densities}.
\newblock {\em Ann. Statist.}, 29(5):1233--1263.

\bibitem[Ghosal and van~der Vaart, 2007]{ghosal2007posterior}
Ghosal, S. and van~der Vaart, A.~W. (2007).
\newblock {Posterior convergence rates of Dirichlet mixtures at smooth
  densities}.
\newblock {\em Ann. Statist.}, 35(2):697--723.

\bibitem[Gin{\'e} and Nickl, 2016]{gine2016mathematical}
Gin{\'e}, E. and Nickl, R. (2016).
\newblock {\em Mathematical Foundations of Infinite-Dimensional Statistical
  Models}.
\newblock Cambridge University Press.

\bibitem[Glorot et~al., 2011]{glorot2011deep}
Glorot, X., Bordes, A., and Bengio, Y. (2011).
\newblock Deep sparse rectifier neural networks.
\newblock In {\em Proc. International Conference on Artificial Intelligence and
  Statistics}, pages 315--323.

\bibitem[Goodfellow et~al., 2014]{goodfellow2014generative}
Goodfellow, I., Pouget-Abadie, J., Mirza, M., Xu, B., Warde-Farley, D., Ozair,
  S., Courville, A., and Bengio, Y. (2014).
\newblock Generative adversarial nets.
\newblock In {\em Proc. Neural Information Processing Systems}, pages
  2672--2680.

\bibitem[Hastie et~al., 2009]{hastie2009elements}
Hastie, T., Tibshirani, R., and Friedman, J.~H. (2009).
\newblock {\em The Elements of Statistical Learning: Data Mining, Inference,
  and Prediction}.
\newblock Springer, New York.

\bibitem[Ho et~al., 2020]{ho2020denoising}
Ho, J., Jain, A., and Abbeel, P. (2020).
\newblock Denoising diffusion probabilistic models.
\newblock In {\em Proc. Neural Information Processing Systems}, volume~33,
  pages 6840--6851.

\bibitem[Huang et~al., 2021]{huang2021error}
Huang, J., Jiao, Y., Li, Z., Liu, S., Wang, Y., and Yang, Y. (2021).
\newblock An error analysis of generative adversarial networks for learning
  distributions.
\newblock {\em ArXiv:2105.13010}.

\bibitem[Hyv{\"a}rinen, 2005]{hyvarinen2005estimation}
Hyv{\"a}rinen, A. (2005).
\newblock Estimation of non-normalized statistical models by score matching.
\newblock {\em J. Mach. Learn. Res.}, 6:695--708.

\bibitem[Imaizumi and Fukumizu, 2022]{imaizumi2022adaptive}
Imaizumi, M. and Fukumizu, K. (2022).
\newblock Adaptive approximation and generalization of deep neural network with
  intrinsic dimensionality.
\newblock {\em J. Mach. Learn. Res.}, 23(111):1--54.

\bibitem[Jordan et~al., 1999]{jordan1999introduction}
Jordan, M.~I., Ghahramani, Z., Jaakkola, T.~S., and Saul, L.~K. (1999).
\newblock An introduction to variational methods for graphical models.
\newblock {\em Mach. Learn.}, 37(2):183--233.

\bibitem[Kim et~al., 2020]{kim2020casting}
Kim, D., Hwang, J., and Kim, Y. (2020).
\newblock {On casting importance weighted autoencoder to an EM algorithm to
  learn deep generative models}.
\newblock In {\em Proc. International Conference on Artificial Intelligence and
  Statistics}, pages 2153--2163. PMLR.

\bibitem[Kingma and Ba, 2015]{kingma2014adam}
Kingma, D.~P. and Ba, J. (2015).
\newblock Adam: A method for stochastic optimization.
\newblock In {\em Proc. International Conference on Learning Representations}.

\bibitem[Kingma and Welling, 2014]{kingma2013auto}
Kingma, D.~P. and Welling, M. (2014).
\newblock {Auto-encoding variational Bayes}.
\newblock In {\em Proc. International Conference on Learning Representations},
  pages 1--14.

\bibitem[Kruijer et~al., 2010]{kruijer2010adaptive}
Kruijer, W., Rousseau, J., and van~der Vaart, A. (2010).
\newblock {Adaptive Bayesian density estimation with location-scale mixtures}.
\newblock {\em Electron. J. Stat.}, 4:1225--1257.

\bibitem[Kundu and Dunson, 2014]{kundu2014latent}
Kundu, S. and Dunson, D.~B. (2014).
\newblock Latent factor models for density estimation.
\newblock {\em Biometrika}, 101(3):641--654.

\bibitem[Li et~al., 2017]{li2017mmd}
Li, C.-L., Chang, W.-C., Cheng, Y., Yang, Y., and P{\'o}czos, B. (2017).
\newblock {MMD GAN: Towards deeper understanding of moment matching network}.
\newblock In {\em Proc. Neural Information Processing Systems}, pages
  2203--2213.

\bibitem[Liang, 2021]{liang2021well}
Liang, T. (2021).
\newblock How well generative adversarial networks learn distributions.
\newblock {\em J. Mach. Learn. Res.}, 22(228):1--41.

\bibitem[Luise et~al., 2020]{luise2020generalization}
Luise, G., Pontil, M., and Ciliberto, C. (2020).
\newblock {Generalization properties of optimal transport GANs with latent
  distribution learning}.
\newblock {\em ArXiv:2007.14641}.

\bibitem[Maas et~al., 2013]{maas2013rectifier}
Maas, A.~L., Hannun, A.~Y., Ng, A.~Y., et~al. (2013).
\newblock Rectifier nonlinearities improve neural network acoustic models.
\newblock {\em ICML Workshop on Deep Learning for Audio, Speech, and Language
  Processing}.

\bibitem[Mhaskar, 1996]{mhaskar1996neural}
Mhaskar, H.~N. (1996).
\newblock Neural networks for optimal approximation of smooth and analytic
  functions.
\newblock {\em Neural Comput.}, 8(1):164--177.

\bibitem[Mroueh et~al., 2018]{mroueh2018sobolev}
Mroueh, Y., Li, C.-L., Sercu, T., Raj, A., and Cheng, Y. (2018).
\newblock Sobolev gan.
\newblock In {\em Proc. International Conference on Learning Representations},
  pages 1--27.

\bibitem[M{\"u}ller, 1997]{muller1997integral}
M{\"u}ller, A. (1997).
\newblock Integral probability metrics and their generating classes of
  functions.
\newblock {\em Adv. in Appl. Probab.}, 29(2):429--443.

\bibitem[Neal, 2011]{neal2011mcmc}
Neal, R.~M. (2011).
\newblock {MCMC using Hamiltonian dynamics}.
\newblock In Brooks, S., Gelman, A., Jones, G., and Meng, X.-L., editors, {\em
  Handbook of Markov Chain Monte Carlo}, chapter~5, pages 113--162. CRC Press.

\bibitem[Ohn and Kim, 2019]{ohn2019smooth}
Ohn, I. and Kim, Y. (2019).
\newblock Smooth function approximation by deep neural networks with general
  activation functions.
\newblock {\em Entropy}, 21(7):627.

\bibitem[Oko et~al., 2023]{oko2023diffusion}
Oko, K., Akiyama, S., and Suzuki, T. (2023).
\newblock Diffusion models are minimax optimal distribution estimators.
\newblock {\em ArXiv:2303.01861}.

\bibitem[Pati et~al., 2011]{pati2011posterior}
Pati, D., Bhattacharya, A., and Dunson, D.~B. (2011).
\newblock Posterior convergence rates in non-linear latent variable models.
\newblock {\em ArXiv:1109.5000}.

\bibitem[Pedregosa et~al., 2011]{pedregosa2011scikit}
Pedregosa, F., Varoquaux, G., Gramfort, A., Michel, V., Thirion, B., Grisel,
  O., Blondel, M., Prettenhofer, P., Weiss, R., Dubourg, V., et~al. (2011).
\newblock Scikit-learn: Machine learning in python.
\newblock {\em J. Mach. Learn. Res.}, 12:2825--2830.

\bibitem[Ramachandran et~al., 2017]{ramachandran2017searching}
Ramachandran, P., Zoph, B., and Le, Q.~V. (2017).
\newblock Searching for activation functions.
\newblock {\em ArXiv:1710.05941}.

\bibitem[Rezende and Mohamed, 2015]{rezende2015variational}
Rezende, D. and Mohamed, S. (2015).
\newblock Variational inference with normalizing flows.
\newblock In {\em Proc. International Conference on Machine Learning}, pages
  1530--1538. PMLR.

\bibitem[Rezende et~al., 2014]{rezende2014stochastic}
Rezende, D.~J., Mohamed, S., and Wierstra, D. (2014).
\newblock Stochastic backpropagation and approximate inference in deep
  generative models.
\newblock In {\em Proc. International Conference on Machine Learning}, pages
  1278--1286.

\bibitem[Schmidt-Hieber, 2020]{schmidt2020nonparametric}
Schmidt-Hieber, J. (2020).
\newblock {Nonparametric regression using deep neural networks with ReLU
  activation function}.
\newblock {\em Ann. Statist.}, 48(4):1875--1897.

\bibitem[Schreuder et~al., 2021]{schreuder2021statistical}
Schreuder, N., Brunel, V.-E., and Dalalyan, A. (2021).
\newblock Statistical guarantees for generative models without domination.
\newblock In {\em Proc. Algorithmic Learning Theory}, pages 1051--1071. PMLR.

\bibitem[Shen et~al., 2013]{shen2013adaptive}
Shen, W., Tokdar, S.~T., and Ghosal, S. (2013).
\newblock {Adaptive Bayesian multivariate density estimation with Dirichlet
  mixtures}.
\newblock {\em Biometrika}, 100(3):623--640.

\bibitem[Singh et~al., 2018]{singh2018nonparametric}
Singh, S., Uppal, A., Li, B., Li, C.-L., Zaheer, M., and P{\'o}czos, B. (2018).
\newblock Nonparametric density estimation with adversarial losses.
\newblock In {\em Proc. Neural Information Processing Systems}, pages
  10246--10257.

\bibitem[Sohl-Dickstein et~al., 2015]{sohl2015deep}
Sohl-Dickstein, J., Weiss, E., Maheswaranathan, N., and Ganguli, S. (2015).
\newblock Deep unsupervised learning using nonequilibrium thermodynamics.
\newblock In {\em Proc. International Conference on Machine Learning}, pages
  2256--2265.

\bibitem[Song and Ermon, 2019]{song2019generative}
Song, Y. and Ermon, S. (2019).
\newblock Generative modeling by estimating gradients of the data distribution.
\newblock {\em Proc. Neural Information Processing Systems}, 32:1--13.

\bibitem[Song et~al., 2020]{song2020sliced}
Song, Y., Garg, S., Shi, J., and Ermon, S. (2020).
\newblock Sliced score matching: A scalable approach to density and score
  estimation.
\newblock In {\em Uncertainty in Artificial Intelligence}, pages 574--584.
  PMLR.

\bibitem[Song et~al., 2021]{song2021scorebased}
Song, Y., Sohl-Dickstein, J., Kingma, D.~P., Kumar, A., Ermon, S., and Poole,
  B. (2021).
\newblock Score-based generative modeling through stochastic differential
  equations.
\newblock In {\em Proc. International Conference on Learning Representations},
  pages 1--36.

\bibitem[Tang and Yang, 2023]{tang2023minimax}
Tang, R. and Yang, Y. (2023).
\newblock Minimax rate of distribution estimation on unknown submanifold under
  adversarial losses.
\newblock {\em Ann. Statist.}, 51(3):1282 -- 1308.

\bibitem[Tang and Yang, 2024]{tang2024adaptivity}
Tang, R. and Yang, Y. (2024).
\newblock Adaptivity of diffusion models to manifold structures.
\newblock In {\em Proc. International Conference on Artificial Intelligence and
  Statistics}.

\bibitem[Telgarsky, 2016]{telgarsky2016benefits}
Telgarsky, M. (2016).
\newblock Benefits of depth in neural networks.
\newblock In {\em Proc. Conference on Learning Theory}, pages 1517--1539.

\bibitem[Tsybakov, 2008]{tsybakov2008introduction}
Tsybakov, A.~B. (2008).
\newblock {\em Introduction to Nonparametric Estimation}.
\newblock Springer, New York.

\bibitem[Uppal et~al., 2019]{uppal2019nonparametric}
Uppal, A., Singh, S., and P{\'o}czos, B. (2019).
\newblock {Nonparametric density estimation and convergence of GANs under Besov
  IPM losses}.
\newblock In {\em Proc. Neural Information Processing Systems}, pages
  9089--9100.

\bibitem[van~der Vaart and Wellner, 1996]{van1996weak}
van~der Vaart, A.~W. and Wellner, J.~A. (1996).
\newblock {\em Weak Convergence and Empirical Processes}.
\newblock Springer.

\bibitem[Villani, 2008]{villani2008optimal}
Villani, C. (2008).
\newblock {\em Optimal Transport: Old and New}.
\newblock Springer.

\bibitem[Vincent, 2011]{vincent2011connection}
Vincent, P. (2011).
\newblock A connection between score matching and denoising autoencoders.
\newblock {\em Neural Comput.}, 23(7):1661--1674.

\bibitem[Virtanen et~al., 2020]{virtanen2020scipy}
Virtanen, P., Gommers, R., Oliphant, T.~E., Haberland, M., Reddy, T.,
  Cournapeau, D., Burovski, E., Peterson, P., Weckesser, W., Bright, J., {van
  der Walt}, S.~J., Brett, M., Wilson, J., Millman, K.~J., Mayorov, N., Nelson,
  A. R.~J., Jones, E., Kern, R., Larson, E., Carey, C.~J., Polat, {\.I}., Feng,
  Y., Moore, E.~W., {VanderPlas}, J., Laxalde, D., Perktold, J., Cimrman, R.,
  Henriksen, I., Quintero, E.~A., Harris, C.~R., Archibald, A.~M., Ribeiro,
  A.~H., Pedregosa, F., {van Mulbregt}, P., and {SciPy 1.0 Contributors}
  (2020).
\newblock {SciPy} 1.0: Fundamental algorithms for scientific computing in
  python.
\newblock {\em Nature Methods}, 17:261--272.

\bibitem[Wong and Shen, 1995]{wong1995probability}
Wong, W.~H. and Shen, X. (1995).
\newblock {Probability inequalities for likelihood ratios and convergence rates
  of sieve MLEs}.
\newblock {\em Ann. Statist.}, 23(2):339--362.

\bibitem[Yang and Zhou, 2023]{yang2023optimal}
Yang, Y. and Zhou, D.-X. (2023).
\newblock Optimal rates of approximation by shallow ${{\rm ReLU}}^k$ neural
  networks and applications to nonparametric regression.
\newblock {\em ArXiv:2304.01561}.

\bibitem[Yarotsky, 2017]{yarotsky2017error}
Yarotsky, D. (2017).
\newblock {Error bounds for approximations with deep ReLU networks}.
\newblock {\em Neural Networks}, 94:103--114.

\end{thebibliography}

		\clearpage
		\section*{Checklist}

		\begin{enumerate}

			\item For all models and algorithms presented, check if you include:
			\begin{enumerate}
				\item A clear description of the mathematical setting, assumptions, algorithm, and/or model. [Yes, Section \ref{sec:model}, \ref{sec:main} and \ref{sec:alternative}]
				\item An analysis of the properties and complexity (time, space, sample size) of any algorithm. [Yes, Section \ref{sec:experiments}]
				\item (Optional) Anonymized source code, with specification of all dependencies, including external libraries. [No]
			\end{enumerate}

			\item For any theoretical claim, check if you include:
			\begin{enumerate}
				\item Statements of the full set of assumptions of all theoretical results. [Yes, Section \ref{sec:main} and \ref{sec:alternative}]
				\item Complete proofs of all theoretical results. [Yes, Supplements]
				\item Clear explanations of any assumptions. [Yes, Section \ref{sec:main} and \ref{sec:alternative}]     
			\end{enumerate}

			\item For all figures and tables that present empirical results, check if you include:
			\begin{enumerate}
				\item The code, data, and instructions needed to reproduce the main experimental results (either in the supplemental material or as a URL). [Yes, Section \ref{sec:experiments}]
				\item All the training details (e.g., data splits, hyperparameters, how they were chosen). [Yes, Section \ref{sec:experiments}]
				\item A clear definition of the specific measure or statistics and error bars (e.g., with respect to the random seed after running experiments multiple times).  [Yes, Section \ref{sec:experiments}]
				\item A description of the computing infrastructure used. (e.g., type of GPUs, internal cluster, or cloud provider). [No]
			\end{enumerate}
			
			\item If you are using existing assets (e.g., code, data, models) or curating/releasing new assets, check if you include:
			\begin{enumerate}
				\item Citations of the creator If your work uses existing assets. [Not Applicable]
				\item The license information of the assets, if applicable. [Not Applicable]
				\item New assets either in the supplemental material or as a URL, if applicable. [Not Applicable]
				\item Information about consent from data providers/curators. [Not Applicable]
				\item Discussion of sensible content if applicable, e.g., personally identifiable information or offensive content. [Not Applicable]
			\end{enumerate}
			
			\item If you used crowdsourcing or conducted research with human subjects, check if you include:
			\begin{enumerate}
				\item The full text of instructions given to participants and screenshots. [Not Applicable]
				\item Descriptions of potential participant risks, with links to Institutional Review Board (IRB) approvals if applicable. [Not Applicable]
				\item The estimated hourly wage paid to participants and the total amount spent on participant compensation. [Not Applicable]
			\end{enumerate}
			
		\end{enumerate}
		
		\appendix
		\onecolumn
		\section{PROOF OF THEOREM \ref{thm:main} } \label{sec:supp1}

		We first state and prove several lemmas needed for proving Theorem \ref{thm:main}.
		
		\medskip
		
		\begin{Lem}
			\label{sec:discrete approximation} Let $p_0 \in {\cC}^{\beta,L,\tau_0}(\bbR^d)$ be a probability density function satisfying assumptions (Tail 1) and (Tail 2).
			Then, there exist positive constants $C_{1} = {C_{1}}({{\rm all}}), C_{2}={C_{2}}({{\rm all}}), C_{3} = C_{3}({\rm all}), C_{4}={C_{4}}({\rm all})$ and a probability measure $H_0$ supported within $[-a_{\sigma},a_{\sigma}]^d$ such that $d_H(p_0,\phi_{\sigma}*H_{0}) \leq C_{2} \sigma^{\beta}$ and $1-P_0([-a_{\sigma},a_{\sigma}]^{d}) \leq C_3 \sigma^{4\beta+8}$ for every $\sigma \in (0, \min(C_1,1))$, where $a_{\sigma} = C_4 \left\{ \log(1/\sigma) \right\}^{\tau_3}$.
		\end{Lem}
		\begin{proof}
			This is a re-statement of Lemma 9.11 in \cite{ghosal2017fundamentals} except for the assertion $1-P_0([-a_{\sigma}, a_{\sigma}]^{d}) \leq C_3 \sigma^{4\beta+8}$, which can be easily derived from the proof of Lemma 9.11.
		\end{proof}
		
		\medskip
		
		\begin{Lem}
			\label{sec:mixture hellinger}
			For any functions $\bff , \bg : [0,1]^{d_0} \rightarrow \bbR^d$ and $\sigma > 0$,
			\bean
			d_{H}^2(p_{\bff,\sigma},p_{\bg,\sigma}) \leq \frac{\| \bff - \bg \|_2^2}{8\sigma^2}.
			\eean
		\end{Lem}
		
		\begin{proof}
			Note that $p_{\bff,\sigma}(\bx) = \int_{[0,1]^{d_0}} \phi_{\sigma}(\bx - \bff(\bz))\d \bz$ and $p_{\bg,\sigma}(\bx) = \int_{[0,1]^{d_0}} \phi_{\sigma}(\bx - \bg(\bz))\d \bz$. We can rewrite squared Hellinger distance as
			\bean
			&& d_{H}^2(p_{\bff,\sigma}, p_{\bg,\sigma}) 
			\\
			&&= \int \left\{ p_{\bff,\sigma}(\bx) +  p_{\bg,\sigma}(\bx) - 2\sqrt{p_{\bff,\sigma}(\bx)}\sqrt{  p_{\bg,\sigma}(\bx)} \right\} \d \bx
			\\
			&&= \int \left[ \int_{[0,1]^{d_0}}  \left\{ \phi_{\sigma}(\bx - \bff(\bz)) + \phi_{\sigma}(\bx - \bg(\bz)) \middle\} \d \bz - 2\sqrt{p_{\bff,\sigma}(\bx)}\sqrt{  p_{\bg,\sigma}(\bx) } \right. \right] \d \bx.
			\eean
			\Holder's inequality implies that 
			\bean
			\int_{[0,1]^{d_0}} \sqrt{\phi_{\sigma} (\bx - \bff(\bz))} \sqrt{\phi_{\sigma} (\bx - \bg(\bz))} \d \bz
			\leq
			\sqrt{p_{\bff,\sigma}(\bx)}\sqrt{ p_{\bg,\sigma}(\bx)}.
			\eean
			Hence,
			\bean
			&& d_{H}^2(p_{\bff,\sigma}, p_{\bg,\sigma}) 
			\\
			&&\leq \int \int_{[0,1]^{d_0}} \left\{ \phi_{\sigma}(\bx - \bff(\bz)) + \phi_{\sigma}(\bx - \bg(\bz)) - 2\sqrt{\phi_{\sigma} (\bx - \bff(\bz))} \sqrt{\phi_{\sigma} (\bx - \bg(\bz))} \right\} \d \bz \d \bx
			\\
			&&= \int \int_{[0,1]^{d_0}} \left\{ \sqrt{\phi_{\sigma} (\bx - \bff(\bz))} - \sqrt{\phi_{\sigma} (\bx - \bg(\bz))}   \right\}^2 \d \bz \d \bx
			\\
			&&= \int_{[0,1]^{d_0}} d_{H}^2 \big( \phi_\sigma(\cdot-\bff(\bz)), \phi_\sigma(\cdot-\bg(\bz)) \big) \d \bz,
			\eean
			where the last equality holds by Fubini's theorem. The squared Hellinger distance between $\cN(\mu_1,\Sigma_1)$ and $\cN(\mu_2,\Sigma_2)$ is known as 
			\bean
			1-\frac{\operatorname{det}\left(\Sigma_1\right)^{1 / 4} \operatorname{det}\left(\Sigma_2\right)^{1 / 4}}{\operatorname{det}\left(\frac{\Sigma_1+\Sigma_2}{2}\right)^{1 / 2}} \exp \left\{-\frac{1}{8}\left(\mu_1-\mu_2\right)^T\left(\frac{\Sigma_1+\Sigma_2}{2}\right)^{-1}\left(\mu_1-\mu_2\right)\right\}.
			\eean
			Using that, we have
			\bean
			d_{H}^2 (p_{\bff,\sigma}, p_{\bg,\sigma}) &\leq& \int_{[0,1]^{d_0}}
			\left\{ 1 -  \ \exp \left( - \frac{\| \bff(\bz) - \bg(\bz) \|_2^2}{8\sigma^2} \right) \right\} \d \bz
			\\
			&\leq&  \frac{\| \bff - \bg \|_2^2}{8\sigma^2}
			\eean
			since $1-\exp(-x) \leq x$ for all $x \in \bbR$.
		\end{proof}
		
		\medskip
		
		\begin{Cor}
			\label{sec:finite mixture} Define $p : \bbR^d \rightarrow \bbR$ as $p(\bx) = \sum_{i=1}^{n} w^{(i)} \phi_{\sigma}(\bx-\bx^{(i)})$ with $\sum_{i=1}^n w^{(i)} = 1, w^{(i)} > 0, $ and $\bx^{(i)} \in \bbR^d$ for each $i$. For $0 < w' \leq w^{(1)}$ and $\bx' \in \bbR^d$, define  $p' : \bbR^d \rightarrow \bbR$ as $p'(\bx) =  w' \phi_{\sigma}(\bx-\bx')  + (w^{(1)}-w') \phi_{\sigma}(\bx-\bx^{(1)}) + \sum_{i=2}^{n} w^{(i)} \phi_{\sigma}(\bx-\bx^{(i)})$. Then, 
			\bean
			d_H^2(p,p') \leq  w' \frac{ \|\bx' - \bx^{(1)}\|_2^2}{8\sigma^2} \leq w'\frac{ \|\bx'\|_2^2  + \|\bx^{(1)}\|_2^2}{4\sigma^2}.
			\eean
		\end{Cor}
		\begin{proof}
			Let $q^{(0)} = 0$ and $q^{(i)} = q^{(i-1)} + w^{(i)}$ for $i \in \{1,\ldots,n\}$.
			Consider functions $\bg,\bg' : [0,1] \rightarrow \bbR^d$ such that $\bg(0)=\bx^{(1)}, \bg'(0) = \bx',$
			\bean
			\bg(z) &=& \sum_{i=1}^n \bx^{(i)} 1_{(q^{(i-1)},q^{(i)}]}(z)  \quad {{\rm and}}
			\\
			\bg'(z) &=& \bx' 1_{ (q^{(0)},w' ]}(z) 
			+ \bx^{(1)} 1_{(w',q^{(1)} ]}(z) +
			\sum_{i=2}^n \bx^{(i)} 1_{(q^{(i-1)},q^{(i)}]}(z)  
			. 
			\eean
			Then, $p_{\bg,\sigma} = p$ and $p_{\bg',\sigma} = p'$. By Lemma \ref{sec:mixture hellinger},
			\bean
			d_H^2(p,p') \leq \frac{\| \bg - \bg' \|_2^2}{8\sigma^2} = w' \frac{ \|\bx' - \bx^{(1)}\|_2^2}{8\sigma^2}.
			\eean
			Since $\|\bx' - \bx^{(1)}\|_2^2 \leq  2\|\bx'\|_2^2  + 2\|\bx^{(1)}\|_2^2$, we obtain the results.
		\end{proof}
		
		\medskip

		The proof of Lemma \ref{sec:bracket} below is almost the same as that of Lemma 1 in \cite{chae2023likelihood}, which is limited to a fixed $F$.
		
		\begin{Lem}
			\label{sec:bracket} Suppose that $0 < \sigma_{\min} \leq 1/\sqrt{2}, \sigma_{\max} \geq 1$ and $F \geq 1$. Let $\cG$ be a class of functions from $[0,1]^{d_0}$ to $\bbR^d$ such that $\|\bg\|_{\infty} \leq F$ for every $\bg \in \cG$. Let $\cP = \{p_{\bg,
				\sigma} : \bg \in \cG, \sigma \in [\sigma_{\min},\sigma_{\max}] \}$. Then, there exist positive constants $C_{5} = C_{5}(d), C_{6} = C_{6}(d)$ and $C_{7} = C_{7}(d)$ such that
			\bean
			\log N_{[]}(\delta,\cP,d_H) &\leq& \log N\left( \frac{C_{5} \delta^4\sigma_{\min}^{d+2} }{F \sigma_{\max}^{2d} \left[ \left\{ \log (\sigma_{\max}/\sigma_{\min}) \right\}^d + F^{2d} \right] }, \cG, \|\cdot\|_{\infty} \right) 
			\\
			&+& \log \left( \frac{ C_{6} \sigma_{\max}^{2d+1} \left[ \left\{ \log (\sigma_{\max}/\sigma_{\min}) \right\}^d + F^{2d} \right]}{ \delta^4\sigma_{\min}^{d+1} } \right)
			\eean
			for $0 < \delta \leq C_7$.
		\end{Lem}
		\begin{proof}
			For $\bg_1,\bg_2 \in \cG$ and $\sigma \in [\sigma_{\min},\sigma_{\max}]$ with $\|\bg_1 - \bg_2\|_{\infty} \leq \eta_1$, we have
			\bean
			&& p_{\bg_1,\sigma}(\bx) - p_{\bg_2,\sigma}(\bx)
			\\
			&&= \int_{[0,1]^{{d_0}}} \phi_{\sigma}(\bx - \bg_1(\bz)) \left\{1- \frac{\phi_{\sigma}(\bx-\bg_2(\bz))}{\phi_{\sigma}(\bx-\bg_1(\bz))}  \right\} \d \bz
			\\
			&&= \int_{[0,1]^{d_0}} \phi_{\sigma}(\bx - \bg_1(\bz)) \left\{1- \exp \middle(-\frac{\| \bx - \bg_2(\bz)\|_2^2 - \|\bx - \bg_1(\bz)\|_2^2}{2\sigma^2} \middle)  \right\} \d \bz
			\\
			&&\leq \int_{[0,1]^{d_0}} \phi_{\sigma}(\bx - \bg_1(\bz)) \left\{ \frac{\| \bx - \bg_2(\bz)\|_2^2 - \|\bx - \bg_1(\bz)\|_2^2}{2\sigma^2}  \right\}  \d \bz
			\\
			&&= \int_{[0,1]^{d_0}} \phi_{\sigma}(\bx - \bg_1(\bz)) \left[ \frac{\|\bg_2(\bz)\|_2^2 - \|\bg_1(\bz)\|_2^2 - 2\bx^{{\rm T}}\left\{\bg_2(\bz)-\bg_1(\bz)\right\}}{2\sigma^2} \right]    \d \bz.
			\eean
			For $\bg_1(\bz) = \left(\left\{\bg_1(\bz)\right\}_{1},\ldots,\left\{\bg_1(\bz)\right\}_{d} \right)$ and $\bg_2(\bz) = \left(\left\{\bg_2(\bz)\right\}_{1},\ldots,\left\{\bg_1(z)\right\}_{d} \right)$, note that
			$
			\|\bg_2(\bz)\|_2^2 - \|\bg_1(\bz)\|_2^2 = \sum_{i=1}^{d}  \left\{\bg_2(\bz)\right\}_{i}^2 - \left\{\bg_2(\bz)\right\}_{i}^2 \leq 2  F d \eta_1.
			$
			Also, it holds that $ \vert \bx^{{\rm T}} (\bg_2(\bz)-\bg_1(\bz)) \vert \leq \|\bx \|_2 \| \bg_2(\bz) - \bg_1(\bz) \|_2  \leq \|\bx\|_2 \sqrt{d} \eta_1 $. Simple calculation yields that $\|\bx\|_2 \leq \| \bx - \bg_1(\bz) \|_2 + \| \bg_1(\bz) \|_2 \leq 1+\|\bx-\bg_1(\bz)\|_2^2 +  F\sqrt{d}$.
			Combining with the last display, we have
			\bean
			&& p_{\bg_1,\sigma}(\bx) - p_{\bg_2,\sigma}(\bx) 
			\\
			&& \leq \int_{[0,1]^{d_0}} \phi_{\sigma}(\bx - \bg_1(\bz)) \left( \frac{ 2Fd\eta_1 + 2\sqrt{d}\eta_1+ 2\sqrt{d} \eta_1 \| \bx - \bg_1(\bz) \|_2^2 + 2F d \eta_1}{2\sigma^2} \right) \d \bz
			\\
			&&\leq \eta_1 \int_{[0,1]^{d_0}} \phi_{\sigma}(\bx - \bg_1(\bz)) \left(\frac{ 2 F d  +\sqrt{d} }{\sigma^2}
			+\frac{\sqrt{d} \| \bx - \bg_1(\bz) \|_2^2}{\sigma^2}  \right) \d \bz
			\\
			&&= \eta_1 (2\pi \sigma^2)^{-d/2} \int_{[0,1]^{d_0}}  \exp\left(- \frac{\|\bx-\bg_1(\bz)\|_2^2}{2\sigma^2} \right) \left(\frac{ 2 F d  +\sqrt{d} }{\sigma^2}
			+\frac{\sqrt{d} \| \bx - \bg_1(\bz) \|_2^2}{\sigma^2}  \right) \d \bz
			\\
			&&\leq \eta_1 (2\pi \sigma^2)^{-d/2} \left(\frac{2Fd+\sqrt{d}}{\sigma^2} + \frac{2\sqrt{d}}{e} \right),
			\eean
			where the last inequality holds because for any $t > 0, e^{-t} \leq 1$ and $ te^{-t} \leq 1/e$.
			Then, there exists a positive constant $D_{1} = D_{1}(d)$ where the last display is further bounded by $\eta_1 D_{1} \sigma_{\min}^{-d-2}F$ for every $F \geq 1$ and $0 < \sigma_{\min} \leq \sqrt{e}$.
			\\
			Also, for $\sigma_1,\sigma_2 \in [\sigma_{\min},\sigma_{\max}]$ and $\bg \in \cG$ with $|\sigma_1-\sigma_2| \leq \eta_2$, we have
			\bean
			&& p_{\bg,\sigma_1}(\bx) - p_{\bg,\sigma_2}(\bx)
			\\
			&&= \int_{[0,1]^{d_0}} \phi_{\sigma_1}(\bx - \bg(\bz)) \left[1-  \exp \left\{- \frac{\|\bx-\bg(\bz)\|_2^2}{2} \left(\frac{1}{\sigma_2^2} - \frac{1}{\sigma_1^2} \right) 
			+d \log \frac{\sigma_1}{\sigma_2} \right\}    \right] \d \bz
			\\
			&&\leq \int_{[0,1]^{d_0}} \phi_{\sigma_1}(\bx - \bg(\bz)) \left\{ \frac{\|\bx-\bg(\bz)\|_2^2}{2} \left(\frac{1}{\sigma_2^2} - \frac{1}{\sigma_1^2} \right) 
			- d \log \frac{\sigma_1}{\sigma_2}    \right\} \d \bz.
			\eean
			Simple calculation yields that $|\sigma_2^{-2} - \sigma_{1}^{-2}| \leq \sigma_1^{-2} \sigma_2^{-2} (\sigma_1+\sigma_2) \eta_2$ and $| \log(\sigma_2/\sigma_1)| \leq \eta_2 / \min(\sigma_1,\sigma_2)$.
			Combining with the last display, we have
			\bean
			&& p_{\bg,\sigma_1}(\bx) - p_{\bg,\sigma_2}(\bx)
			\\
			&&\leq \eta_2 \int_{[0,1]^{d_0}} \phi_{\sigma_1}(\bx - \bg(\bz)) \left\{ \frac{(\sigma_1+\sigma_2)  \|\bx-\bg(\bz)\|_2^2}{2\sigma_1^2 \sigma_2^2}  
			+ \frac{d}{\min(\sigma_1,\sigma_2)}   \right\} \d \bz
			\\
			&&= \eta_2 (2\pi\sigma_1^2)^{-d/2} \int_{[0,1]^{d_0}} \exp\left(- \frac{\|\bx-\bg(\bz)\|_2^2}{2\sigma_1^2} \right) \left\{ \frac{(\sigma_1+\sigma_2)  \|\bx-\bg(\bz)\|_2^2}{2\sigma_1^2 \sigma_2^2}  
			+ \frac{d}{\min(\sigma_1,\sigma_2)}   \right \}\d \bz
			\\
			&&\leq  \eta_2 (2\pi\sigma_1^2)^{-d/2} \left\{ \frac{\sigma_1+\sigma_2}{e \sigma_2^2}  
			+ \frac{d}{\min(\sigma_1,\sigma_2)}   \right\}.
			\eean
			Then, there exist a positive constant $D_{2} = D_{2}(d)$ where the last display is further bounded by $\eta_2 D_{2} \sigma_{\min}^{-d-1}$.
			\\
			Given $0 < \epsilon < 1$, set $\eta_1 = \epsilon /(4 D_{1} \sigma_{\min}^{-d-2}F)$ and  $\eta_2 = \epsilon /(4 D_{2} \sigma_{\min}^{-d-1})$.
			Suppose $\{\bg_1,\ldots,\bg_{N_1} \}$ be $\eta_1$-covering set of $\cG$ and $\{ \sigma_1,\ldots,\sigma_{N_{2}} \}$ be $\eta_2$-covering set of $[\sigma_{\min},\sigma_{\max}]$.
			Then, $\{p_{\bg_i,\sigma_j}: i \in \{1,\ldots,N_1\}, j \in \{1,\ldots,N_{2}\} \}$ forms an $\epsilon/2$-covering set of $(\cP,\|\cdot\|_{\infty})$ for every $F \geq 1$ and $0< \sigma_{\min} \leq 1 \leq \sigma_{\max}$.
			Define $l_{ij}$ and $u_{ij}$ as
			\bean
			l_{ij}(\bx) = \max \left\{ p_{\bg_i,\sigma_j}(\bx)-\epsilon/2,0 \right\} \quad {{\rm and}} \quad u_{ij}(\bx) = \min\left\{p_{\bg_i,\sigma_j}(\bx)+\epsilon/2,H(\bx) \right\}
			\eean
			for each $(i,j)$, where $H(\bx) = \sup_{p \in \cP } p(\bx)$. 
			Note that
			\bean
			H(\bx) &\leq& (2\pi\sigma_{\min}^2)^{-d/2} \sup_{\|\by\|_{\infty} \leq F} \exp \left( -\frac{\|\bx-\by\|_2^2}{2\sigma_{\max}^2} \right)
			\\\
			&\leq& (2\pi \sigma_{\min}^2)^{-d/2} \exp \left(- \frac{\|\bx\|_2^2 - 2dF^2}{4\sigma_{\max}^2} \right) 
			= \left(\frac{\sqrt{2}\sigma_{\max}}{\sigma_{\min}}\right)^d \exp \left(\frac{dF^2}{2\sigma_{\max}^2} \right) \phi_{\sqrt{2}\sigma_{\max}}(\bx),
			\eean
			where the last inequality holds because $\|\bx-\by\|_2^2 \geq \| \bx \|_2^2 /2 - \|\by\|_2^2 \geq \|\bx\|_2^2 /2 - d F^2$ for $\|\by\|_{\infty} \leq F$.
			For any $t > 0$, Gaussin tail bound implies that $\int_{\|\bx\|_{\infty} > t} \phi_{\sqrt{2}\sigma_{\max}}(\bx)\d \bx \leq d e^{-t^2/(4\sigma_{\max}^2)}$.
			Since $\sigma_{\max} \geq 1$, we have that 
			\bean
			\int_{\| \bx \|_{\infty} > B} H(\bx) \d \bx &\leq& 
			\left(\frac{\sqrt{2}\sigma_{\max}}{\sigma_{\min}}\right)^d \exp \left(\frac{dF^2}{2} \right) \int_{\|\bx\|_{\infty} > B} \phi_{\sqrt{2}\sigma_{\max}} (\bx) \d \bx
			\\
			&\leq& d \left(\frac{\sqrt{2}\sigma_{\max}}{\sigma_{\min}}\right)^d \exp \left(\frac{dF^2}{2} -\frac{B^2}{4\sigma_{\max}^2} \right) = \epsilon
			\eean
			where 
			\bean
			B = 2\sigma_{\max} \left(\log \frac{1}{\epsilon} + d \log \frac{\sigma_{\max}}{\sigma_{\min}} + \frac{d}{2}\log 2+ \frac{dF^2}{2} + \log d \right)^{1/2}.
			\eean
			Hence,
			\bean
			\| u_{ij} - l_{ij} \|_1 &=& \int_{\bbR^d} \left\{ u_{ij}(\bx)-l_{ij}(\bx) \right\} \d \bx
			\\
			&\leq& \int_{\|\bx\|_{\infty} \leq B} \epsilon \ \d \bx + \int_{\|\bx\|_{\infty} > B} H(\bx) \d \bx \leq \left\{(2B)^{d}+1 \right\}\epsilon.
			\eean
			Define $\delta = \sqrt{\epsilon \{(2B)^d+1 \}}$. Since $d_H^2(u_{ij},l_{ij}) \leq \| u_{ij} - l_{ij} \|_1$, we have
			\bean
			N_{[]}(\delta,\cP,d_H) \leq N_{[]}(\delta^2,\cP,\| \cdot \|_1) \leq N_1 N_{2} \leq \frac{\sigma_{\max}-\sigma_{\min}}{\eta_2} N(\eta_1,\cG,\|\cdot\|_{\infty})
			\eean
			for every $F \geq 1$ and $0< \sigma_{\min} \leq 1 \leq \sigma_{\max}$.
			\\
			There exists a positive constant $D_{3} = D_{3}(d)$ such that for $0 < \sigma_{\min} \leq 1/\sqrt{2}$ and $ 1 \leq \sigma_{\max}$,
			\bean
			\delta^2 = \epsilon (2^d B^d + 1) &\leq& \epsilon D_{3} \sigma_{\max}^d \left[  \left\{\log \left(\frac{1}{\epsilon} \right) \right\}^{d/2} + \left\{\log \left(\frac{\sigma_{\max}}{\sigma_{\min}}\right)\right\}^{d/2} + F^d \right]
			\eean
			Since $\epsilon \leq \epsilon \{ \log (1/\epsilon) \}^{d/2} \leq \sqrt{\epsilon}$ for $\epsilon < \epsilon_1$, where $\epsilon_1 = \epsilon_1(d) < 1$ is a constant, we have
			\bean
			\delta^2 \leq \sqrt{\epsilon} D_{3} \sigma_{\max}^d \left[ \left\{ \log \left(\frac{\sigma_{\max}}{\sigma_{\min} }\right) \right\}^{d/2} + F^d \right].
			\eean
			Hence,
			\bean
			\eta_1 = \epsilon \frac{ \sigma_{\min}^{d+2}}{4 D_{1} F} &\geq& \frac{ \sigma_{\min}^{d+2} \delta^4}{4 D_{1} F D_{3}^2 \sigma_{\max}^{2d} \left[ \left\{ \log (\sigma_{\max}/\sigma_{\min}) \right\}^{d/2} + F^d \right]^2} \\
			&\geq& \frac{D_{4} \delta^4\sigma_{\min}^{d+2} }{F \sigma_{\max}^{2d} \left[ \left\{ \log (\sigma_{\max}/\sigma_{\min}) \right\}^{d} + F^{2d} \right]}
			\eean
			and
			\bean
			\eta_2 = \epsilon \frac{ \sigma_{\min}^{d+1}}{4 D_{2}} &\geq& 
			\frac{ \sigma_{\min}^{d+1} \delta^4}{4 D_{2} D_{3}^2 \sigma_{\max}^{2d} \left[ \left\{ \log (\sigma_{\max}/\sigma_{\min}) \right\}^{d/2} + F^d \right]^2}
			\\ &\geq& \frac{D_{5} \delta^4\sigma_{\min}^{d+1} }{ \sigma_{\max}^{2d} \left[ \left\{ \log (\sigma_{\max}/\sigma_{\min}) \right\}^{d} + F^{2d} \right]}
			\eean
			for $D_{4} = D_{4}(d), D_{5} = D_{5}(d)$ and $0< \delta \leq D_{6} < 1$, where $D_{6} = D_{6}(d)$.
			The assertion follows by re-defining constants.
		\end{proof}

		\medskip
		The proof of Lemma \ref{sec:nnentropy} below is a straightforward extension of Lemma 5 in \cite{schmidt2020nonparametric}, which can only be applied to the case of $M=1$.
		
		\begin{Lem}
			\label{sec:nnentropy}
			For any $\delta > 0$ and $\bd = (1,d_1,d)$, we have
			\bean
			\log N(\delta,\cG(1,\infty,\bd,M),\| \cdot\|_{\infty}) \leq d_1(d+2)\log\left(\frac{ 8M^2d_1}{\delta} \right).
			\eean
		\end{Lem}
		\begin{proof}
			Let $z \mapsto \bg(z) = W_2\rho_{\bb}W_1 z$ and $z \mapsto \bg'(z) =  W'_2\rho_{\bb'}W'_1 z$ with $\bg,\bg' \in \cG(1,\infty,\bd,M)$. Given $\epsilon > 0$, assume that all parameter values of $\bg$ and $\bg'$ are at most $\epsilon$ away from each other. 
			Then,
			\bean
			\| \bg(z) - \bg'(z) \|_{\infty}
			&\leq&
			\| W_2\rho_{\bb}W_1 z - W_2\rho_{\bb'}W_1'  z\|_{\infty} + 
			\| W_2\rho_{\bb'}W'_1 z - W'_2\rho_{\bb'}W'_1  z \|_{\infty}
			\\
			&=& \| W_2 (\rho_{\bb}W_1 z -\rho_{\bb'}W_1'  z)\|_{\infty} + 
			\| (W_2-W'_2) \rho_{\bb'}W'_1 z  \|_{\infty}
			\\
			&\leq& Md_1 \| \rho_{\bb}W_1 z - \rho_{\bb'}W'_1  z\|_{\infty} + \epsilon d_1 \| \rho_{\bb'}W'_1 z  \|_{\infty},
			\eean
			where the last inequlaity holds because for any matrix $W \in  \bbR^{d \times d_1}$ and $\bx \in \bbR^{d_1}$, we have $\| W \bx \|_{\infty} \leq d_1 \|W\|_{\infty}\| \bx\|_{\infty} $.
			It holds that $\| \rho_{\bb}W_1 z - \rho_{\bb'}W'_1  z\|_{\infty} \leq \|(W_1-W'_1)z \|_{\infty} + \|\bb - \bb'\|_{\infty}$ and $\| \rho_{\bb'}W'_1 z  \|_{\infty} \leq \|W_1' z\|_{\infty} + \|\bb'\|_{\infty}$. Combining with the last diplay, we have
			\bean
			\| \bg(z) - \bg'(z) \|_{\infty} 
			&\leq&  M d_1 \left\{ \|(W_1-W'_1)z \|_{\infty} + \|\bb - \bb'\|_{\infty} \right\} + \epsilon d_1 \left\{ \|W_1' z\|_{\infty} + \|\bb'\|_{\infty} \right\}
			\\
			&\leq& \epsilon M d_1 \left( \|z\|_{\infty} + 1 \right) + \epsilon M d_1 \left( \|z\|_{\infty} + 1 \right)
			\\
			&=& 2\epsilon M d_1 \left( \|z\|_{\infty} + 1 \right)
			\\
			&\leq&  4 \epsilon M d_1 .
			\eean
			Note that the total number of parameters in $\bg$ is equal to $d_1(d+2)$. Define $\delta = 4 \epsilon M d_1 $. Then,
			\bean
			N\left(\delta,\cG(1,\infty,\bd,M),\| \cdot\|_{\infty} \right)
			&\leq& N \left(\epsilon,\left[-M,M\right]^{d_1(d+2)}, \| \cdot \|_{\infty} \right)
			\\
			&\leq&  \left( \frac{2M}{\epsilon} \right)^{d_1(d+2)} = \left( \frac{8M^2 d_1}{\delta} \right)^{d_1(d+2)}.
			\eean
			The assertion follows by taking a logarithm.
		\end{proof}

		\bigskip
		
		\textit{Proof of Theorem \ref{thm:main}.} 
		We will apply Theorem 4 of \cite{wong1995probability} with $\alpha = 0+$.
		Let $c_{1},\ldots,c_{4}$ be the same positive constants defined in Theorem 1 of \cite{wong1995probability}. These constants can be chosen, for example, as $c_{1} = 1/24, c_2=2/26001, c_3 = 10$ and $ c_4 = (2/3)^{5/2}/512$.
		By Theorem 4 of \cite{wong1995probability}, it suffices to prove that
		\bean
		\int_{\epsilon_n^2 / 2^8}^{\sqrt{2}\epsilon_n} \sqrt{\log N_{[]}(\delta/c_{3},\cP,d_H )} \ \d \delta \leq c_{4} \sqrt{n} \epsilon_n^2
		\eean
		and there exist $\bg_* \in \cG(1,F,\bd,M)$ and $\sigma_* \in [\sigma_{\min}, \sigma_{\max}]$ satisfying
		\begin{align}
			\int \log\left(\frac{p_0(\bx)}{p_{\bg_*, \sigma_*}(\bx)} \right) \d P_0(\bx) &< \frac{1}{4}c_{1} \epsilon_n^2 \label{eqn:approx1}
			\\
			\int \left\{ \log \left(\frac{p_0(\bx)}{p_{\bg_*,\sigma_*}(\bx)} \right) \right\}^2 \d P_0(\bx) &< \frac{1}{4}c_{1} \epsilon_n^2\log n, \label{eqn:approx2}
		\end{align}
		for every $n \geq \widetilde C_1$, where $\widetilde C_1$ and $\widetilde C_2$ are large enough constants and $\epsilon_n$ is defined as in Theorem \ref{thm:main}.

		To derive \eqref{eqn:approx1} and \eqref{eqn:approx2}, we firstly approximate $p_0$ by Gaussian mixture densities and then construct ReLU networks to approximate the mixing measure.
		Techniques approximating $p_0$ by Gaussian mixtures are originally developed by \cite{shen2013adaptive} and slightly refined in \cite{ghosal2017fundamentals}.
		
		Let $C_1,\ldots,C_4$ be constants in Lemma \ref{sec:discrete approximation} and $a_\sigma = C_4\left\{\log(1/\sigma)\right\}^{\tau_3}$.
		Let $\sigma \in [\sigma_{\min}, \sigma_{\max}]$ be small enough as described below.
		By Lemma \ref{sec:discrete approximation}, if $\sigma \leq \min(C_1,1)$, there exists a probability measure $H_0$ supported within $[-a_{\sigma}, a_{\sigma}]^{d}$ such that 
		\be \label{eq:hellinger1}
		d_H(p_0,\phi_{\sigma} * H_{0}) \leq C_{4} \sigma^{\beta}. 
		\ee 
		If, furthermore, $\sigma$ is small enough so that $a_\sigma/\sigma \geq 1$, then Lemma 9.12 of \cite{ghosal2017fundamentals} implies that there exist positive constants $ D_{1} = D_{1}(d,\beta), D_{2} = D_{2}(d,\beta)$ and discrete probability measure $\widetilde{H}_{0}(\cdot) = \sum_{t=1}^{N_0} w^{(t)} \delta_{\bx^{(t)}}(\cdot)$, where $\delta_\bx(\cdot)$ denotes the Dirac measure at $\bx$, supported inside $[-a_{\sigma},a_{\sigma}]^d$ such that
		\bean
		N_{0} \leq D_{1} a_{\sigma}^d \sigma^{-d} \left\{ \log (1/\sigma) \right\}^d = D_{1} C_{4}^d \sigma^{-d}  \left\{\log(1/\sigma)\right\}^{\tau_3 d+d}
		\eean
		and
		\be \label{eq:hellinger2}
		d_H(\phi_{\sigma} * H_{0} , \phi_{\sigma} * \widetilde{H}_{0}) \leq \| \phi_{\sigma} * H_{0} - \phi_{\sigma} * \widetilde{H}_{0} \|_1^{1/2} \leq D_{2} \sigma^{\beta} \left\{\log(1/\sigma)\right\}^{d/4}.
		\ee
		Moreover, $\widetilde H_0$ can be constructed so that $\bx^{(1)}, \ldots, \bx^{(N_0)}$ are distinct, $w^{(t)} > 0$ and
		\be \label{eq:grid}
		\left\{\bx^{(1)},\ldots,\bx^{(N_0)} \right\}
		\subseteq
		\left\{(n_1,\ldots,n_d)\sigma^{2\beta+1} : n_i \in \bbZ, i=1, \ldots, d \right\} \cap \left[-a_{\sigma},a_{\sigma} \right]^d.
		\ee
		Without loss of generality, we may assume that $w^{(1)} \geq \cdots \geq w^{(N_0)}$. Let 
		\bean
		N_1 = \left|\left\{t : w^{(t)} \geq \sigma^{2\beta+2d+2}, ~ t =1, \ldots, N_0 \right\}\right|,
		\eean
		where $|\cdot|$ denotes the cardinality.
		If $\sigma$ is small enough, we have $1 \leq N_1 \leq N_0$. 
		Let $\widetilde{H}_1(\cdot) = w^{(+)} \delta_{\bx^{(1)}}(\cdot) +  \sum_{t=2}^{N_1} w^{(t)} \delta_{\bx^{(t)}} (\cdot)$ where $w^{(+)} = w^{(1)} + \sum_{t > N_1} w^{(t)}
		$.
		Corollary \ref{sec:finite mixture} implies that 
		\be \begin{split} \label{eq:hellinger3}
			&d_H(\phi_{\sigma} * \widetilde{H}_0, \phi_{\sigma} * \widetilde{H}_1) 
			\leq \sqrt{\frac{d}{2}} (N_0 - N_1) a_\sigma \sigma^{\beta+d}
			\\
			&\leq \sqrt{\frac{d}{2}} a_{\sigma} N_0 \sigma^{\beta+d} 
			\leq \sqrt{\frac{d}{2}} D_1 C_{4}^{d+1} \sigma^{\beta} \left\{\log(1/\sigma)\right\}^{\tau_3 d+\tau_3+d}.
		\end{split} \ee
		
		For $t =1, \ldots, N_1$, let $U_t$ be the intersection of $[-a_\sigma, a_\sigma]^d$ and $\bx^{(t)} + B(\sigma^{2\beta+1}/3)$, the $\ell_2$-ball with the radius $\sigma^{2\beta+1}/3$ centered on $\bx^{(t)}$.
		Since $\bx^{(1)}, \ldots, \bx^{(N_1)}$ are on grids \eqref{eq:grid}, $U_1, \ldots, U_{N_1}$ are mutually disjoint.
		One can extend $\{U_1, \ldots, U_{N_1}\}$ to $\{U_1, \ldots, U_{N_2}\}$ so that the latter forms a partition of $[-a_\sigma, a_\sigma]^d$ and the $\ell_2$-diameter of $U_t$ is at most $\sigma$ for all $t \leq N_2$.
		Since 
		\bean
		N\left(\frac{\sigma}{2}, [-a_\sigma, a_\sigma]^d, \|\cdot\|_2\right) 
		\leq N\left(\frac{\sigma}{2\sqrt{d}}, [-a_\sigma, a_\sigma]^d, \|\cdot\|_\infty \right)
		\leq \left( \frac{4\sqrt{d}a_{\sigma}}{\sigma} \right)^{d},
		\eean
		one may construct a partition so that $N_2 \leq N_1 + (4\sqrt{d}a_{\sigma}/\sigma)^{d}$.
		Hence,
		\bean
		N_{2} \leq N_1 + (4\sqrt{d}a_{\sigma}/\sigma)^{d} \leq  N_0 + (4\sqrt{d}a_{\sigma}/\sigma)^{d}  \leq D_{3} \sigma^{-d} \left\{\log(1/\sigma) \right\}^{\tau_3d+d},
		\eean
		where $D_{3} = D_{3}(C_{4},D_{1})$.
		Define $\widetilde{\bx}^{(t)}$ as $\widetilde{\bx}^{(t)} = \bx^{(t)}$ for $t \leq N_1$ and choose $\widetilde{\bx}^{(t)} \in U_t$ for $N_1 < t \leq N_2$.
		Let $\widetilde{H}_2(\cdot) = \sum_{t=1}^{N_{2}} \widetilde{w}^{(t)} \delta_{\widetilde{\bx}^{(t)}}(\cdot)$
		where $\widetilde{w}^{(1)} = w^{(+)} - (N_{2}-N_1)\sigma^{2\beta+2d+2}, \widetilde{w}^{(t)} = w^{(t)}$ for $1 < t \leq N_1$ and $\widetilde{w}^{(t)} = \sigma^{2\beta+2d+2}$ for $N_1 < t \leq N_2$.
		Since $w^{(+)} \geq w^{(1)} \geq N_0^{-1} \gtrsim \sigma^d \{\log (1/\sigma)\}^{-\tau_3d - d}$ and $N_2 \lesssim \sigma^{-d} \left\{\log(1/\sigma) \right\}^{\tau_3 d + d}$, we have $\widetilde{w}^{(1)} \geq \sigma^{2\beta+2d+2}$ for small enough $\sigma$.
		Corollary \ref{sec:finite mixture} implies that
		\be \begin{split} \label{eq:hellinger4}
			&d_H(\phi_{\sigma} * \widetilde{H}_1, \phi_{\sigma} * \widetilde{H}_2) \leq \sqrt{\frac{d}{2}} (N_2 - N_1) a_\sigma \sigma^{\beta+d} 
			\\
			&< \sqrt{\frac{d}{2}} a_{\sigma} N_2 \sigma^{\beta+d} 
			\leq  \sqrt{\frac{d}{2}} D_{3} C_4 \sigma^{\beta} \left\{\log(1/\sigma) \right\}^{\tau_3 d+\tau_3+d}.
		\end{split} \ee
		
\begin{figure}
			\centering
			\includegraphics[width = 0.8\textwidth]{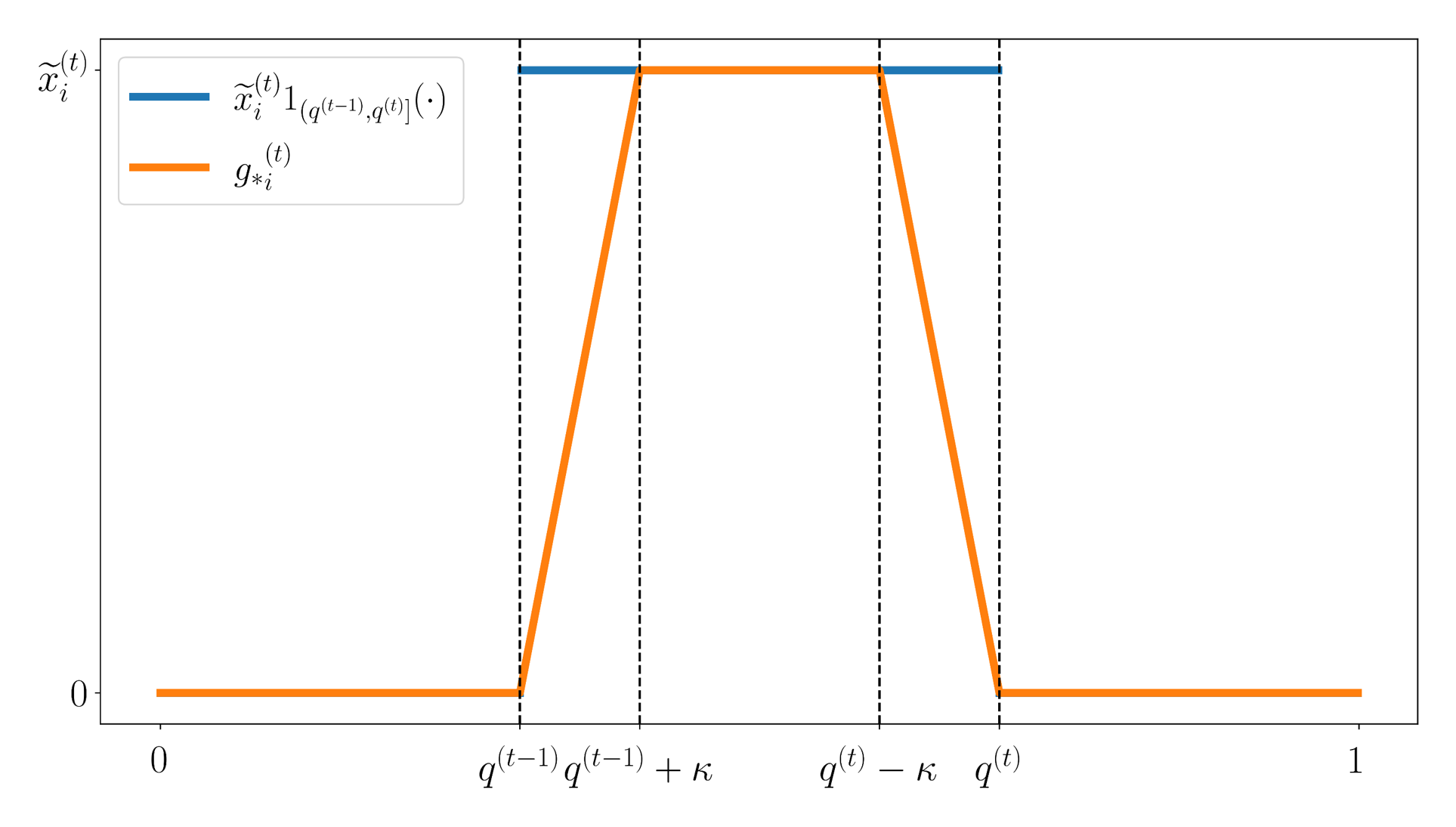}
			\caption{Step function approximation with ReLU network}
			\label{fig:relu}
		\end{figure}

		Consider function $\widetilde{\bg} = (\widetilde{g}_{1},\ldots,\widetilde{g}_d) : [0,1] \rightarrow [-a_{\sigma},a_{\sigma}]^{d}$ such that
		\bean
		\widetilde{\bg}(0) = \widetilde{\bx}^{(1)} \quad {\rm and} \quad \widetilde{\bg}(z) = \sum_{t=1}^{N_{2}} \widetilde{\bx}^{(t)} 1_{(q^{(t-1)},q^{(t)}]}(z),
		\eean
		where $q^{(0)} = 0$ and $q^{(t)} = \sum_{s=1}^{t} \widetilde{w}^{(s)}$ for $t \leq N_2$. Then, $\phi_{\sigma}*\widetilde{H}_2 = p_{\widetilde{\bg},\sigma}$.
		Note that $\widetilde{g}_i : [0,1] \rightarrow [-a_{\sigma},a_{\sigma}]$ is a step function with $\widetilde{g}_i(0) = \widetilde{x}_i^{(1)}$ and $\widetilde{g}_i(z) = \sum_{t=1}^{N_{2}} \widetilde{x}_i^{(t)} 1_{(q^{(t-1)},q^{(t)} ]}(z)$.
		Define ${g_*}_i^{(t)} \in \cG(1,a_{\sigma},(1,4,1),\kappa^{-1})$ as
		\be \begin{split} \label{eq:nn_aprox}
			{g_*}_i^{(t)}(z) &= \widetilde{x}_i^{(t)} \left[\max\left\{0,\frac{ 1}{\kappa}\left(z - q^{(t-1)}\right) \right\}
			- \max\left\{0,\frac{ 1}{\kappa}\left(z - \left(q^{(t-1)} + \kappa \right) \right) \right\} \right.
			\\
			&\  \left. - \max\left\{0,\frac{1}{\kappa}\left(z - \left(q^{(t)} - \kappa \right) \right) \right\}
			+ \max\left\{0,\frac{1}{\kappa}\left(z - q^{(t)}\right) \right\} \right],
		\end{split} \ee 
		which approximates $\widetilde{x}_i^{(t)} 1_{(q^{(t-1)},q^{(t)} ]}(\cdot)$ as described in Figure \ref{fig:relu}, where $\kappa = \sigma^{2\beta+2d+3}/2$ and we have $ \kappa^{-1} \geq a_{\sigma}$ for small enough $\sigma$.
		Define ${g_*}_i \in \cG(1,a_{\sigma},(1,4 N_{2},1),\kappa^{-1})$ as ${g_*}_i(z) =\sum_{t=1}^{N_{2}} {g_*}_i^{(t)}(z)$ for each $i \in \{1,\ldots,d\}$ and define
		\bean
		\bg_* \in \cG\left(a_{\sigma},(1,4N_{2},d),{\kappa}^{-1} \right) \quad {\rm as} \quad \bg_*(z) = \left({g_*}_{1}(z),\ldots,{g_*}_{d}(z) \right).
		\eean
		Then,
		\bean
		\|\widetilde{\bg} - \bg_* \|_2^2 
		&=& \int_{[0,1]} \sum_{i=1}^d \left\{\widetilde{g}_i(z) - {g_*}_i(z) \right\}^2 \d z 
		\\
		&=& \int_{[0,1]} \sum_{i=1}^d \left[\sum_{t=1}^{N_{2}} \left\{ \widetilde{x}_i^{(t)}1_{(q^{(t-1)},q^{(t)}]}(z) - {g_*}_i^{(t)}(z) \right\} \right]^2 \d z 
		\\
		&=& \int_{[0,1]} \sum_{i=1}^d \sum_{t=1}^{N_{2}} 
		\left\{\widetilde{x}_i^{(t)}1_{(q^{(t-1)},q^{(t)}]}(z)  -{g_*}_i^{(t)}(z) \right\}^2 \d z 
		\\
		&=&  2\sum_{i=1}^d \sum_{t=1}^{N_{2}} \int_{[0,\kappa]}
		\left(\widetilde{x}_i^{(t)} -\frac{\widetilde{x}_i^{(t)}}{\kappa}z\right)^2 \d z 
		\\
		&=& \frac{2}{3}\sum_{i=1}^d \sum_{t=1}^{N_{2}} \kappa \left\{ \widetilde{x}_i^{(t)} \right\}^2 \leq  \frac{2}{3}d N_{2} \kappa a_{\sigma}^2.
		\eean
		Combining with Lemma \ref{sec:mixture hellinger}, we have
		\bean
		&& d_H(\phi_{\sigma} * \widetilde{H}_2,p_{\bg_*,\sigma}) = d_H(p_{\widetilde{\bg},\sigma},p_{\bg_*,\sigma})
		\\
		&& \leq \frac{\|\widetilde{\bg} - \bg_* \|_2}{2\sqrt{2}\sigma} 
		=
		\frac{\sqrt{d N_{2} \kappa}a_{\sigma}}{2\sqrt{3}\sigma}
		\leq D_{4} \sigma^{\frac{2\beta + d+1}{2}} \left\{\log(1/\sigma) \right\}^{\frac{\tau_3 d +2\tau_3 +d}{2}},
		\eean
		where $D_{4} = \sqrt{dD_{3}}C_{4}/(2\sqrt{6})$.
		Combining \eqref{eq:hellinger1}, \eqref{eq:hellinger2}, \eqref{eq:hellinger3} and \eqref{eq:hellinger4} with the last display, we have
		\bean
		d_H(p_0,p_{\bg_*,\sigma}) 
		&\leq& C_{4}\sigma^{\beta} + D_{2} \sigma^{\beta} \left\{\log(1/\sigma) \right\}^{\frac{d}{4}} +
		\sqrt{\frac{d}{2}} D_{1} C_{4}^{d+1} \sigma^{\beta} \left\{\log(1/\sigma)\right\}^{\tau_3 d+\tau_3+d} 
		\\ &\ & + \sqrt{\frac{d}{2}} D_{3} C_{4} \sigma^{\beta} \left\{\log(1/\sigma)\right\}^{\tau_3 d+\tau_3+d}
		+ D_{4} \sigma^{\frac{2\beta + d+1}{2}} \left\{\log(1/\sigma) \right\}^{\frac{\tau_3 d +2\tau_3 +d}{2}}
		\\
		&\leq& D_{5} \sigma^{\beta} \left\{\log(1/\sigma)\right\}^{\tau_3 d+\tau_3+d},
		\eean
		where $D_{5} = C_{4}+D_{2} +\sqrt{d/2}D_{1}C_{4}^{d+1}+ \sqrt{d/2}D_{3}C_{4}+D_{4}$.
		
		For any $\bx \in [-a_{\sigma},a_{\sigma}]^{d}$, there exists $s \in \{1,\ldots,N_{2} \}$ such that $\bx \in U_{s}$. Since $\widetilde{\bx}^{(s)} \in U_{s}$ and $\| \widetilde{\bx}^{(s)} - \bx \|_2 \leq \sigma$, we have
		\be \begin{split} \label{eq:D6}
			p_{\bg_*,\sigma}(\bx) &\geq \int_{ \{ \|\bx - \bg_*(z) \|_2 \leq \sigma \} } \phi_{\sigma}(\bx-\bg_*(z)) \d z 
			\\
			&= (2\pi)^{-\frac{d}{2}}\int_{ \{\|\bx - \bg_*(z) \|_2 \leq \sigma \}}  \sigma^{-d} \exp\left(-\frac{ \|\bx-\bg_*(z) \|_2^2}{2\sigma^2} \right)\d z 
			\\
			&\geq (2\pi)^{-\frac{d}{2}} \sigma^{-d} \int_{ \{ \widetilde{\bx}^{(s)} = \bg_*(z) \}} e^{-\frac{1}{2}}  \d z
			= (2\pi)^{-\frac{d}{2}} e^{-\frac{1}{2}} \sigma^{-d} (\widetilde{w}^{(s)}-2\kappa)
			\\
			&> D_{6} \sigma^{2\beta+d+2},
		\end{split}\ee
		where the last inequality holds because $\min_t \widetilde{w}^{(t)} \geq \sigma^{2\beta+2d+2}$ and $D_{6} = (2\pi)^{-d/2}e^{-1/2}/2$ .
		For any $\bx \in \bbR^d$ with $\| \bx \|_{\infty} > a_{\sigma}$, we have
		\bean
		p_{\bg_*,\sigma}(\bx) &=& (2\pi)^{-\frac{d}{2}} \sigma^{-d} \int_{ \| \bg_*(z) \|_{\infty} \leq  a_\sigma}\exp\left(-\frac{ \|\bx-\bg_*(z) \|_2^2}{2\sigma^2} \right)\d z
		\\
		&\geq&  (2\pi)^{-\frac{d}{2}} \sigma^{-d} \exp\left(-\frac{ 2d \|\bx\|_2^2}{\sigma^2}\right),
		\eean
		where the last inequality holds because $\|\bx - \bg_*(z) \|_2^2 \leq 2\|\bx\|_2^2 + 2 \|\bg_*(z)\|_2^2 \leq 2\|\bx\|_2^2 + 2d \|\bg_*(z)\|_{\infty}^2 \leq 4d \|\bx\|_2^2$.
		Combining with (Tail 2) assumption, it follows that $
		p_0(\bx)/p_{\bg_*,\sigma}(\bx) \leq  \tau_1 (2\pi)^{d/2} \sigma^{d} \exp\left(2d\|\bx\|_2^2/\sigma^2 \right). $
		Hence,
		\bean
		\log \left(\frac{p_0(\bx)}{p_{\bg_*,\sigma}(\bx)} \right) \leq D_{7} +  \frac{2d \|\bx\|_2^2}{\sigma^2},
		\eean
		where $D_{7} = \log \tau_1 + d \log (2\pi)/2$.
		
		Assumption (Tail 2) and \eqref{eq:D6} implies that $p_{\bg_*,\sigma}(\bx) / p_{0}(\bx) > \lambda$ for all $\bx \in[-a_\sigma, a_\sigma]^d$, where $\lambda = D_6\sigma^{2\beta+d+2} / \tau_1$. 
		It follows that
		$\{\bx: p_{\bg_*,\sigma}(\bx)/p_0(\bx) \leq \lambda, \bx \in \bbR^{d} \} \subseteq \{\bx: \|\bx\|_{\infty}> a_{\sigma}, \bx \in \bbR^{d} \}$.
		Hence, 
		\bean
		&& \int_{ \left\{ \frac{p_{\bg_*,\sigma}(\bx)}{p_0(\bx)} \leq \lambda \right\}} \left\{\log \left(\frac{p_0(\bx)}{p_{\bg_*,\sigma}(\bx)} \right) \right\}^2 \d P_0(\bx)
		\\
		&& \leq \int_{ \{ \|\bx\|_{ \infty} > a_{\sigma} \}} \left\{ \log
		\left(\frac{p_0(\bx)}{p_{\bg_*,\sigma}(\bx)} \right) \right\}^2 \d P_0(\bx)
		\\
		&& \leq 2 \int_{ \{ \|\bx\|_{\infty} > a_{\sigma} \}} \left( D_{7}^2 + \frac{4d^2 \|\bx\|_2^4}{\sigma^4} \right) \d P_0(\bx)
		\\
		&& \leq 2D_{7}^2 \left\{1 - P_{0}\left( \left[-a_{\sigma},a_{\sigma} \right]^{d} \right) \right\} 
		+ \frac{8d^2}{\sigma^4} \int_{ \{ \|\bx\|_{\infty} > a_{\sigma} \} } \|\bx\|_2^4 \  \d P_0(\bx)
		\\
		&&\leq 2D_{7}^2 \left\{1 - P_{0}\left( \left[-a_{\sigma},a_{\sigma} \right]^{d} \right) \right\}
		+ \frac{8d^2}{\sigma^4} \left\{ \bbE\left[\|\bX\|_2^8 \right] \right\}^{1/2}  \left\{1 - P_{0}\left( \left[-a_{\sigma},a_{\sigma} \right]^{d} \right) \right\}^{1/2} 
		\\
		&& \leq D_{8} \sigma^{2\beta},
		\eean
		where $D_{8} = D_{8}(C_{3},C_{4},D_{7},\tau_1,\tau_2,\tau_3)$ and the last inequality holds by Lemma \ref{sec:discrete approximation} and (Tail 2) assumption.
		Since $\lambda$ is sufficiently small for small enough $\sigma$, Lemma B.2 of \cite{ghosal2017fundamentals} implies that there exist positive constants $D_{9} = D_{9}(D_{6},d,\beta,\tau_1)$ and $D_{10} = D_{10}(D_{5},D_{8},D_{9})$ such that
		\bean
		&&\int \log \left(\frac{p_0(\bx)}{p_{\bg_*,\sigma}(\bx)} \right)^2 \d P_0(\bx)
		\\
		&&\leq d_H^2(p_0,p_{\bg_*,\sigma})\left[12 + 2 \{\log(1/\lambda )\}^2 \right] 
		+ 8\int_{ \left\{\frac{p_{\bg_*,\sigma}(\bx)}{p_0(\bx)} \leq \lambda \right\}}  \left\{\log \left(\frac{p_0(\bx)}{p_{\bg_*,\sigma}(\bx)} \right) \right\}^2 \d P_0(\bx)
		\\
		&&\leq D_{9} d_H^2(p_0,p_{\bg_*,\sigma}) \left\{ \log (1/\sigma) \right\}^2 + 8D_{8}\sigma^{2\beta}
		\\
		&&\leq D_{10} \sigma^{2\beta} \left\{ \log (1/\sigma) \right\}^{2\tau_3d+2\tau_3+2d+2}
		\eean
		and
		\bean
		&& \int \log \left(\frac{p_0(\bx)}{p_{\bg_*,\sigma}(\bx)} \right) \d P_0(\bx) 
		\\ &&\leq
		d_H^2(p_0,p_{\bg_*,\sigma})\left[1 + 2\log(1/\lambda) \right] + 2\int_{ \left\{ \frac{p_{\bg_*,\sigma}(\bx)}{p_0(\bx)} \leq \lambda \right\} }  \log \left(\frac{p_0(\bx)}{p_{\bg_*,\sigma}(\bx)} \right) \d P_0(\bx)
		\\
		&&\leq D_{9} d_H^2(p_0,p_{\bg_*,\sigma}) \log(1/\sigma) + 2 \int_{ \left\{ \frac{p_{\bg_*,\sigma}(\bx)}{p_0(\bx)} \leq \lambda \right\} }  \left\{\log \left(\frac{p_0(\bx)}{p_{\bg_*,\sigma}(\bx)} \right) \right\}^2 \d P_0(\bx)
		\\
		&&\leq D_{9} d_H^2(p_0,p_{\bg_*,\sigma}) \log(1/\sigma) + 2D_{8}\sigma^{2\beta}
		\\
		&&\leq D_{10} \sigma^{2\beta} \left\{ \log(1/\sigma) \right\}^{2\tau_3d+2\tau_3+2d+1}. 
		\eean
		For $\sigma_* \asymp n^{-1/(2\beta+d)}$ with $\sigma_* \in [\sigma_{\min},\sigma_{\max}]$, if $n$ is large enough, we have
		\begin{align}
			\int \log \left(\frac{p_0(\bx)}{p_{\bg_*,{\sigma_*}}(\bx)} \right) \d P_0(\bx)
			&\leq D_{11} n^{-\frac{2\beta}{2\beta+d}} \left(\log n \right)^{2\tau_3d+2\tau_3+2d+1} \quad {\rm and} \label{eq:kl1}
			\\
			\int \log \left(\frac{p_0(\bx)}{p_{\bg_*,\sigma}(\bx)} \right)^2 \d P_0(\bx) &\leq D_{11} n^{-\frac{2\beta}{2\beta+d}} \left(\log n \right)^{2\tau_3d+2\tau_3+2d+2}, \label{eq:kl2}
		\end{align}
		where $D_{11} =D_{10}/(2\beta+d)^{2\tau_3d+2\tau_3+2d+2}$.
		
		Let $C_{5},\ldots,C_{7}$ be constants in Lemma \ref{sec:bracket}.
		Then, there exists a positive constant $D_{12} = D_{12}({\rm all})$ such that for every $\delta \leq C_{7}$ and large enough $n$,
		\bean
		\log N_{[]}(\delta,\cP,d_H) \leq \log N \left(C_{5} D_{12} \delta^4 n^{-\frac{d+3}{2\beta+d}},\cG,\| \cdot\|_{\infty} \right) + \log\left(\frac{C_{6} D_{12} n^{\frac{d+2}{2\beta+d}}}{\delta^4  } \right).
		\eean
		Combining with Lemma \ref{sec:nnentropy}, we have
		\bean
		\log N_{[]}(\delta,\cP,d_H) \leq  D_{13} n^{\frac{d}{2\beta+d}} (\log n)^{\tau_3d+d} \left\{ \log n + \log(1/\delta) \right\},
		\eean
		where $D_{13} = D_{13}({\rm all})$.
		Note that for every $\epsilon \leq \min( c_{3} C_{7} /\sqrt{2}, 1/e)$, we have
		\bean
		&& \int_{\epsilon^2 / 2^8}^{\sqrt{2}\epsilon} \sqrt{\log N_{[]}(\delta /c_{3},\cP,d_H )} \ \d \delta 
		\\
		&&\leq  \sqrt{2}\epsilon \sqrt{D_{13} n^{d/(2\beta+d)} (\log n)^{\tau_3d+d}  \{\log n + \log(c_3 2^8 / \epsilon^2) \}  }
		\\
		&&\leq D_{14}  n^{\frac{d}{4\beta+2d}} (\log n)^{\frac{\tau_3d+d}{2}}  \epsilon   \left\{ \log n + \log (1/\epsilon) \right\}^{\frac{1}{2}},
		\eean
		where $D_{14} = D_{14}(D_{13},d,\beta)$.
		Therefore, for all large enough $n$, the last display holds with $\epsilon = \epsilon_n$ and is further bounded by $c_4 \sqrt{n} \epsilon_n^2$, where
		\bean
		\epsilon_n = D_{15} n^{-\frac{\beta}{2\beta+d}}  \left( \log n \right)^{\frac{2\tau_3 d+2\tau_3 +2d + 1}{2}}
		\eean
		and $D_{15} = D_{15}(D_{11},D_{14},d,\beta,\tau_3)$ is a large enough constant.
		If $D_{15}$ is chosen so that $D_{15} > 4 D_{11}/ c_{1}$, \eqref{eq:kl1} and \eqref{eq:kl2} is further bounded by $c_1 \epsilon_n^2 /4$ and $c_1 \epsilon_n^2 \log n /4$, respectively.
		By re-defining constants, the proof is complete.
		\qed
		
		\newpage
		\section{PROOF OF THEOREM \ref{thm:structured}}
		
		In addition to Section \ref{sec:supp1}, we state and prove lemmas needed for proving Theorem \ref{thm:structured}.
		
		\medskip
		
		\begin{Lem}
			\label{sec:beta2 approximation} Let $p_0 \in {\cC}^{\beta,L,\tau_0}(\bbR^d)$ be a probability density function satisfying an assumption (Tail 1).
			Then, there exist a positive constant $C_{8} = {C_{8}}(d,\beta,L)$ such that $d_H(p_0,\phi_{\sigma}*P_{0}) \leq C_{8} \sigma^{\min(\beta,2)}$ for every $\sigma \in (0, \min( 1/\sqrt{4\tau_0},1)]$.
		\end{Lem}
		\begin{proof}
			For any $\bx,\by \in \bbR^{d}$, Taylor's theorem for multivariate functions yields that
			\bean
			p_0(\bx-\by) - p_0(\bx) = \sum_{1 \leq k. \leq \lfloor \beta \rfloor}  \frac{  (-\by)^{\bk} (D^{\bk}p_0)(\bx)}{\bk!} + \sum_{ k. = \lfloor \beta \rfloor} \frac{(-\by)^{\bk} \left\{ (D^{\bk}p_0)(\bx- t\by) - (D^{\bk}p_0)(\bx) \right\} }{\bk!}
			\eean
			for a suitable $t \in [0,1]$, where $(-\by)^{\bk} = \prod_{i=1}^{d} (-y_i)^{k_i} $ and $\bk! = \prod_{i=1}^d k_i !$.
			Let $m_{\bk} = \prod_{i=1}^{d} m_{k_i}$ and $m_{k_i} = \int \phi(y) y^{k_i} \d y$, where $m_{k_i}$ denote the $k_i$-th moment of standard normal distribution on $\bbR$. In particular, $m_{k_i} = 0$ if $k_i$ is an odd number.
			Since $\int \phi_{\sigma}(\by) \by^{\bk} \d \by =  m_{\bk} \sigma^{k.}$ and combining with the last display, we have
			\bean
			&& (\phi_{\sigma}*P_0)(\bx) - p_0(\bx) = \int_{\bbR^d} \phi_{\sigma}(\by) \left\{ p_0(\bx-\by)-p_0(\bx) \right\} \d\by
			\\
			&& =  \sum_{1 \leq k. \leq \lfloor \beta \rfloor} \frac{  (-\sigma)^{k.}m_{\bk}  (D^{\bk}p_0)(\bx)}{\bk!} + \sum_{ k. = \lfloor \beta \rfloor} \int_{\bbR^d} \frac{\phi_{\sigma}(\by)(-\by)^{\bk}  \left\{ (D^{\bk}p_0)(\bx- t\by) - (D^{\bk}p_0)(\bx) \right\}}{\bk!} \d \by
			\\
			&& \leq \sum_{1 \leq k. \leq \lfloor \beta \rfloor} \frac{  \sigma^{k.}m_{\bk} \left\vert(D^{\bk}p_0)(\bx)\right\vert  }{\bk!} + \sum_{ k. = \lfloor \beta \rfloor} \frac{L(\bx)}{\bk !} \int_{\bbR^d} e^{\tau_0 \| \by\|_2^2} \phi_{\sigma}(\by) {| \by |}^{\bk} \| \by \|_2^{\beta-\lfloor \beta \rfloor} \d \by,
			\eean
			where $|\by|^{\bk} = \prod_{i=1}^d |y_i|^{k_i}$ and the last inequality holds by the definition of $\cC^{\beta,L,\tau_0}(\bbR^{d})$.
			Note that $ e^{\tau_0 \| \by\|_2^2} \phi_{\sigma}(\by) \leq 2^{d/2} \phi_{\sqrt{2}\sigma}(\by)$ because $\tau_0 \leq 1/(4\sigma^2)$. 
			Also, it follows that $|\by|^{\bk} \|\by\|_2^{\beta-\lfloor \beta \rfloor} \leq \sum_{i=1}^d |y_i|^{\beta-\lfloor \beta \rfloor} \prod_{j=1}^{d} |y_j|^{k_j} $ and $\int \phi_{\sqrt{2}\sigma}(\by) |\by|^{\bk} \|\by\|_2^{\beta-\lfloor \beta \rfloor} \d \by \leq (\sqrt{2} \sigma)^{\beta} d (m_{2 \lfloor \beta \rfloor +2})^{d}$ for $k. = \lfloor \beta \rfloor$.
			Since $m_{\bk} = 0$ for $k.=1$ and combining with the last display, it follows that
			\bean
			\left\vert (\phi_{\sigma}*P_0)(\bx) - p_0(\bx) \right\vert \leq D_1 \sigma^{\min(\beta,2)} \left( \sum_{2 \leq k. \leq \lfloor \beta \rfloor} \left\vert(D^{\bk} p_0)(\bx)\right\vert + L(\bx) \right),
			\eean
			where $D_1 = D_1(d,\beta)$.
			Hence,
			\bean
			&& d_H^2(p_0,\phi_{\sigma}*P_0) = \int_{\bbR^d} \left\{\frac{p_0(\bx) - (\phi_{\sigma}*P_0)(\bx) }{ \sqrt{p_0(\bx)} + \sqrt{(\phi_{\sigma}*P_0)(\bx)} } \right\}^2 \d \bx
			\\
			&& \leq D_1^2 \sigma^{2 \min(\beta,2)} \int_{\bbR^d} \left\{ \left( \frac{L(\bx)}{p_0(\bx)} \right)^2 + \sum_{2 \leq k.  \leq \lfloor \beta \rfloor} \left( \frac{|(D^{\bk} p_0)(\bx) |}{p_0(\bx)} \right)^2 \right\} \d P_0(\bx)
			\\
			&& \leq D_2 \sigma^{2 \min(\beta,2)},
			\eean
			where $D_2 = D_2(D_1,\beta,L)$ and the last inequality holds by the assumption (Tail 1).
			The assertion follows by re-defining constants.
		\end{proof}
		\bigskip
		
		\textit{Proof of Theorem \ref{thm:structured}.}
		The proof follows a similar approach to that of Theorem \ref{thm:main}. We will apply Theorem 4 of \cite{wong1995probability} with $\alpha = 0+$.
		Let $c_{1},\ldots,c_{4}$ be the same positive constants defined in Theorem 1 of \cite{wong1995probability}. These constants can be chosen, for example, as $c_{1} = 1/24, c_2=2/26001, c_3 = 10$ and $ c_4 = (2/3)^{5/2}/512$.
		By Theorem 4 of \cite{wong1995probability}, it suffices to prove that
		\bean
		\int_{\epsilon_n^2 / 2^8}^{\sqrt{2}\epsilon_n} \sqrt{\log N_{[]}(\delta/c_{3},\cP,d_H )} \ \d \delta \leq c_{4} \sqrt{n} \epsilon_n^2
		\eean
		and there exist $\bg_* \in \cG(L,F,\bd,M,s)$ and $\sigma_* \in [\sigma_{\min}, \sigma_{\max}]$ satisfying
		\bean
		\int \log\left(\frac{p_0(\bx)}{p_{\bg_*, \sigma_*}(\bx)} \right) \d P_0(\bx) &<& \frac{1}{4}c_{1} \epsilon_n^2
		\\
		\int \left\{ \log \left(\frac{p_0(\bx)}{p_{\bg_*,\sigma_*}(\bx)} \right) \right\}^2 \d P_0(\bx) &<& \frac{1}{4}c_{1} \epsilon_n^2,
		\eean
		for every $n \geq \widetilde C_1$, where $\widetilde C_1$ and $\widetilde C_2$ are large enough constants and $\epsilon_n$ is defined as in Theorem \ref{thm:structured}. 
		
		Let $C_{8}$ be a constant in Lemma \ref{sec:beta2 approximation} and $\sigma \in [\sigma_{\min},\sigma_{\max}]$ be small enough as described below. Combining Lemma \ref{sec:beta2 approximation} and Lemma \ref{sec:mixture hellinger}, if $\sigma \leq \min(1/\sqrt{4\tau_0},1)$, we have
		\begin{align} 
			d_H(p_0,p_{\bg,\sigma}) \leq d_H(p_0,\phi_{\sigma} * P_0) + d_H(p_{\bg_0,\sigma}, p_{\bg,\sigma}) \leq C_8 \sigma^{\widetilde{\beta}} + \frac{\|\bg_0-\bg\|_2}{2\sqrt{2}\sigma} \label{eq:alter_hel1}
		\end{align}
		for any function $\bg : [0,1]^{d} \rightarrow \bbR^{d}$, where the first inequality holds because $\phi_{\sigma}*P_0 = p_{\bg_0,\sigma}$. Lemma 5 of \cite{chae2023likelihood} implies that there exist a constant $D_1 = D_1(\beta,q,\bv,\bt,\bbeta,\tau_6)$ satisfying
		$ \|\bg_0-\bg_{*}\|_{\infty} \leq \sigma^{\widetilde{\beta}+1} $ for some $ \bg_{*} \in \cG\left(\widetilde{L},\infty,\widetilde{\bd},1,\widetilde{s} \right)$,
		where $\widetilde{L} = \lfloor D_1 \log (1/\sigma) \rfloor, \widetilde{\bd} = (d,\widetilde{d}_1,\ldots,\widetilde{d}_1,d)$ with $\widetilde{d}_1 = \lfloor D_1 \sigma^{-\frac{(\widetilde{\beta}+1) t_*}{\beta_*}} \rfloor$ and $\widetilde{s} = D_1 \sigma^{-\frac{(\widetilde{\beta}+1) t_*}{\beta_*}} \log (1/\sigma) $.
		Since $\|\bg_*\|_{\infty} \leq \|\bg_*-\bg_0\|_{\infty} + \| \bg_0\|_{\infty}$ and $\|\bg_0\|_{\infty} \leq \tau_4$ by (Support) assumption, it follows that
		\begin{align}
			\|\bg_0-\bg_{*}\|_{\infty} \leq \sigma^{\widetilde{\beta}+1} \quad \text{for \ some \ } \bg_{*} \in \cG\left(\widetilde{L},\widetilde{F},\widetilde{\bd},1,\widetilde{s} \right), \label{eq:alter_nnaprox}
		\end{align}
		where $\widetilde{F} = \tau_4 + 1$.
		
		Combining \eqref{eq:alter_hel1} and \eqref{eq:alter_nnaprox}, we have
		\begin{align}
			d_H(p_0,p_{\bg_*,\sigma}) \leq (C_8+\sqrt{d/8}) \sigma^{\widetilde{\beta}}, \label{eq:alter_hel2}
		\end{align}
		where the inequality holds because $\|\bg_0-\bg \|_{2} \leq \sqrt{d} \| \bg_0-\bg \|_{\infty}$. 
		
		Assumption (Structured generator) implies that for any $\widetilde{\bx} \in \bbR^{d}$ with $p_0(\widetilde{\bx}) > 0$, there exists $\widetilde{\bz} \in [0,1]^{d}$, such that $\widetilde{\bx} = \bg_0(\widetilde{\bz})$. Note that $\widetilde{\bz}$ does not need to be unique.
		For any $\bz \in [0,1]^{d}$, simple calculation yields that $\| \widetilde{\bx} - \bg_{*}(\bz) \|_2 \leq \| \widetilde{\bx} - \bg_{0}(\bz) \|_2 + \| \bg_{0}(\bz) - \bg_{*}(\bz) \|_2 \leq \sqrt{d} \| \bg_0(\widetilde{\bz}) - \bg_{0}(\bz) \|_{\infty} + \sqrt{d} \|\bg_0-\bg_{*}\|_{\infty}$.
		Combining with \eqref{eq:alter_nnaprox}, we have
		\bean
		\left\{\bz \in [0,1]^{d} : \| \bg_0(\widetilde{\bz}) - \bg_{0}(\bz) \|_{\infty} \leq \sigma \right\}
		&& \subseteq
		\left\{\bz \in [0,1]^{d} : \| \widetilde{\bx} - \bg_{*}(\bz) \|_2 \leq \sqrt{d}(\sigma+ \sigma^{\widetilde{\beta}+1}) \right\}
		\\
		&& \subseteq
		\left\{\bz \in [0,1]^{d} : \| \widetilde{\bx} - \bg_{*}(\bz) \|_2 \leq 2\sqrt{d}\sigma \right\}.
		\eean
		Note that $\bg_0 = \bh_{q} \circ \ldots \circ \bh_{0}$ and $\bh_i = (h_{i1},\ldots,h_{iv_{i+1}})^{\rm T} : [a_i,b_i]^{v_i} \rightarrow [a_i,b_i]^{v_{i+1}}$. Since $h_{ij} \in \cC^{\beta_i}([a_i,b_i]^{t_i};\tau_6)$ and $\|D^{\bk} h_{ij}\|_{\infty} \leq \tau_6$ for $k.=1$, it follows that for any $\bz_1^{(i)},\bz_2^{(i)} \in [a_i,b_i]^{v_i}$, $\|\bh_i(\bz_1^{(i)}) - \bh_i(\bz_2^{(i)}) \|_{\infty} \leq \tau_6 \|\bz_1^{(i)}-\bz_2^{(i)}\|_{\infty}$. Then, simple calculation yields that $\| \bg_0(\bz_1^{(0)}) - \bg_0(\bz_2^{(0)}) \|_{\infty} \leq {\tau_6}^{q+1} \| \bz_1^{(0)} - \bz_2^{(0)} \|_{\infty}$. Combining with the last display, we have
		\bean
		\left\{\bz \in [0,1]^{d} : \| \widetilde{\bz} - \bz \|_{\infty} \leq {\tau_6}^{-(q+1)} \sigma \right\}
		&& \subseteq
		\left\{\bz \in [0,1]^{d} : \| \bg_0(\widetilde{\bz}) - \bg_{0}(\bz) \|_{\infty} \leq \sigma \right\}
		\\
		&& \subseteq
		\left\{\bz \in [0,1]^{d} : \| \widetilde{\bx} - \bg_{*}(\bz) \|_2 \leq 2\sqrt{d}\sigma \right\}.
		\eean
		Hence,
		\bean
		p_{\bg_*,\sigma}(\widetilde \bx) &&\geq \int_{ \{ \|\widetilde \bx - \bg_*(\bz) \|_2 \leq 2\sqrt{d} \sigma \} } \phi_{\sigma}(\widetilde \bx-\bg_*(\bz)) \d \bz 
		\\
		&&= (2\pi)^{-\frac{d}{2}}\int_{ \{\|\widetilde \bx - \bg_*(\bz) \|_2 \leq 2 \sqrt{d} \sigma \}}  \sigma^{-d} \exp\left(-\frac{ \|\widetilde \bx-\bg_*(\bz) \|_2^2}{2\sigma^2} \right)\d \bz 
		\\
		&&\geq (2\pi)^{-\frac{d}{2}} \sigma^{-d}  e^{-2d} \int_{ \{\|\widetilde \bx - \bg_*(\bz) \|_2 \leq 2 \sqrt{d} \sigma \}} \d \bz
		\\
		&&\geq (2\pi)^{-\frac{d}{2}} \sigma^{-d}  e^{-2d} \int_{ \{\|\widetilde \bz - \bz \|_{\infty} \leq {\tau_6}^{-(q+1)} \sigma \}} \d \bz \geq D_2,
		\eean
		where $D_2 = (2\pi)^{-\frac{d}{2}} e^{-2d} {\tau_6}^{-d(q+1)}$ and $\sigma$ is small enough so that ${\tau_6}^{-(q+1)}\sigma \leq 1/2$. Since $\| p_0 \|_{\infty} < \infty$, we have $p_{\bg_*,\sigma}(\bx)/p_0(\bx) > \lambda$ for any $\bx$ with $p_0(\bx) > 0$, where $\lambda = 2^{-1} \min\{ D_2 \| p_0 \|_{\infty}^{-1}, 0.4 \}$. Then, it follows that $\{\bx: p_{\bg_*,\sigma}(\bx) \leq \lambda p_0(\bx) , p_0(\bx) \geq 0 \} \subseteq \{\bx: p_0(\bx) = 0, \bx \in \bbR^d \}$. 
		Lemma B.2 of \cite{ghosal2017fundamentals} and \eqref{eq:alter_hel2} implies that
		\bean
		&&\int \log \left(\frac{p_0(\bx)}{p_{\bg_*,\sigma}(\bx)} \right)^2 \d P_0(\bx)
		\\
		&&\leq d_H^2(p_0,p_{\bg_*,\sigma})\left[12 + 2 \{\log(1/\lambda )\}^2 \right] 
		+ 8\int_{ \left\{\frac{p_{\bg_*,\sigma}(\bx)}{p_0(\bx)} \leq \lambda \right\}}  \left\{\log \left(\frac{p_0(\bx)}{p_{\bg_*,\sigma}(\bx)} \right) \right\}^2 \d P_0(\bx)
		\\
		&&= D_{3}\sigma^{2\widetilde{\beta}}
		\eean
		and
		\bean
		&& \int \log \left(\frac{p_0(\bx)}{p_{\bg_*,\sigma}(\bx)} \right) \d P_0(\bx) 
		\\ &&\leq
		d_H^2(p_0,p_{\bg_*,\sigma})\left[1 + 2\log(1/\lambda) \right] + 2\int_{ \left\{ \frac{p_{\bg_*,\sigma}(\bx)}{p_0(\bx)} \leq \lambda \right\} }  \log \left(\frac{p_0(\bx)}{p_{\bg_*,\sigma}(\bx)} \right) \d P_0(\bx)
		\\
		&&\leq D_{3}\sigma^{2\widetilde{\beta}},
		\eean
		where $D_3 = (C_8+\sqrt{d/8})\left[12 + 2 \{\log(1/\lambda )\}^2 \right]$. 
		
		For $\sigma_* \asymp n^{-\frac{\beta_*}{t_*(\widetilde{\beta}+1) + 2\widetilde{\beta}\beta_* }}$ with $\sigma_* \in [\sigma_{\min}, \sigma_{\max}]$, if $n$ is large enough, we have
		\begin{align}
			\int \log \left(\frac{p_0(\bx)}{p_{\bg_*,{\sigma_*}}(\bx)} \right) \d P_0(\bx)
			&\leq D_{3} n^{-\frac{2\widetilde{\beta}\beta_*}{2\widetilde{\beta}\beta_* + t_*(\widetilde{\beta}+1)}} \quad {\rm and} \label{eq:alter_kl1}
			\\
			\int \log \left(\frac{p_0(\bx)}{p_{\bg_*,\sigma}(\bx)} \right)^2 \d P_0(\bx) &\leq D_{3} n^{-\frac{2\widetilde{\beta}\beta_*}{2\widetilde{\beta}\beta_* + t_*(\widetilde{\beta}+1)}}. \label{eq:alter_kl2}
		\end{align}
		Let $C_{5},\ldots,C_{7}$ be constants in Lemma \ref{sec:bracket}.
		Then, for every $\delta \leq C_{7}$ and large enough $n$ so that $\sigma_{\min}F[\{ \log (1/\sigma_{\min}) \}^{d} + F^{2d}] \leq 1$, we have
		\bean
		\log N_{[]}(\delta,\cP,d_H) \leq \log N \left(C_{5} \delta^4 n^{-\frac{\beta_*(d+3)}{t_*(\widetilde{\beta}+1) + 2\widetilde{\beta}\beta_* }},\cG,\| \cdot\|_{\infty} \right) + \log\left(\frac{C_{6} n^{\frac{\beta_*(d+2)}{t_*(\widetilde{\beta}+1) + 2\widetilde{\beta}\beta_* }} }{\delta^4  } \right).
		\eean
		Lemma 5 of \cite{schmidt2020nonparametric} implies that there exists a constant $D_4 = D_4(d,\beta,\tau_4,q,\bv,\bt,\bbeta,\tau_6)$ such that
		\bean
		\log N_{[]}(\delta,\cP,d_H) \leq  D_{4}n^{\frac{t_*(\widetilde{\beta}+1)}{2\widetilde{\beta} \beta_* + t_*(\widetilde{\beta}+1)}}  \left\{ (\log n)^2 + \log(1/\delta) \right\}.
		\eean
		Note that for every $\epsilon \leq \min( c_{3} C_{7} /\sqrt{2}, 1/e)$, we have
		\bean
		&& \int_{\epsilon^2 / 2^8}^{\sqrt{2}\epsilon} \sqrt{\log N_{[]}(\delta /c_{3},\cP,d_H )} \ \d \delta 
		\\
		&&\leq  \sqrt{2}\epsilon \sqrt{ D_{4}n^{\frac{t_*(\widetilde{\beta}+1)}{2\widetilde{\beta} \beta_* + t_*(\widetilde{\beta}+1)}}  \{(\log n)^2 + \log(c_3 2^8 / \epsilon^2) \}  }
		\\
		&&\leq D_{5}  n^{\frac{t_*(\widetilde{\beta}+1)}{4\widetilde{\beta} \beta_* + 2t_*(\widetilde{\beta}+1)}} \epsilon   \left\{ (\log n)^2 + \log (1/\epsilon) \right\}^{\frac{1}{2}},
		\eean
		where $D_{5} = D_{5}(D_{4},d,\beta,q,\bv,\bt,\bbeta)$.
		Therefore, for all large enough $n$, the last display holds with $\epsilon = \epsilon_n$ and is further bounded by $c_4 \sqrt{n} \epsilon_n^2$, where
		\bean
		\epsilon_n = D_{6} n^{-\frac{\widetilde{\beta}\beta_*}{2\widetilde{\beta}\beta_* + t_*(\widetilde{\beta}+1)}}  (\log n)
		\eean
		and $D_{6} = D_{6}(D_{3},D_{5},d,\beta,q,\bv,\bt,\bbeta)$ is a large enough constant.
		If $D_{6}$ is chosen so that $D_{6} > 4 D_{3}/ c_{1}$, \eqref{eq:alter_kl1} and \eqref{eq:alter_kl2} are further bounded by $c_1 \epsilon_n^2 /4$.
		By re-defining constants, the proof is complete.
		\qed

	\end{document}